\documentclass[9pt]{article}
\usepackage{mathrsfs}
\usepackage{amsthm}
\usepackage{amssymb}
\usepackage{amsmath}
\usepackage{graphicx}
\usepackage{color}
\usepackage{amsfonts}
\usepackage{float}
\usepackage{cite}
\usepackage [latin1]{inputenc}
\usepackage[text={140mm,210mm},left=35mm,vmarginratio=1:1]{geometry}
\newtheorem{theorem}{Theorem}[section]
\newtheorem{remark}{Remark}[section]
\newtheorem{lemma}[theorem]{Lemma}
\newtheorem{proposition}[theorem]{Proposition}

\newtheorem{corollary}[theorem]{Corollary}

\numberwithin{equation}{section}
\normalsize

\begin{document}
\title{\textbf{Scaling limits and fluctuations of a family of $N$-urn branching processes}}

\author{Lirong Ren \thanks{\textbf{E-mail}: 20121634@bjtu.edu.cn \textbf{Address}: School of Mathematics and Statistics, Beijing Jiaotong University, Beijing 100044, China.}\\ Beijing Jiaotong University\\ Xiaofeng Xue \thanks{\textbf{E-mail}: xfxue@bjtu.edu.cn \textbf{Address}: School of Mathematics and Statistics, Beijing Jiaotong University, Beijing 100044, China.}\\ Beijing Jiaotong University}

\date{}
\maketitle

\noindent {\bf Abstract:}
In this paper we are concerned with a family of $N$-urn branching processes, where some particles are put into $N$ urns initially and then each particle gives birth to several new particles in some urn when dies. This model includes the $N$-urn Ehrenfest model and the $N$-urn branching random walk as special cases. We show that the scaling limit of the process is driven by a $C(\mathbb{T})$-valued linear ordinary differential equation and the fluctuation of the process is driven by a generalized Ornstein-Uhlenbeck process in the dual of $C^\infty(\mathbb{T})$, where $\mathbb{T}=(0, 1]$ is the one-dimensional torus. A crucial step for proofs of above main results is to show that numbers of particles in different urns are approximately independent. As applications of our main results, limit theorems of hitting times of the process are also discussed.

\quad


\section{Introduction}\label{section one}
In this paper we are concerned with a family of $N$-urn branching processes. We first introduce the model. Initially, some particles are put into $N$ urns, where $N\geq 1$ is an integer. Let $\mathbb{T}=(0, 1]$ be the one-dimensional torus and $\{\lambda_k\}_{k\geq 0}, \{\psi_k\}_{k\geq 0}$ be positive functions such that $\lambda_k\in C^{\infty}(\mathbb{T}^2)$ and $\psi_k\in C^{\infty}(\mathbb{T})$ for all $k\geq 0$, then our model evolves according to following two rules.

1) For any $1\leq i\neq j\leq N$ and $k\geq 0$, each particle in the $i$th urn dies and meanwhile gives birth to $k$ new particles in the $j$th urn at rate $\frac{1}{N}\lambda_k\left(\frac{i}{N}, \frac{j}{N}\right)$.

2) For any $1\leq i\leq N$ and $k\geq 0$, each particle in the $i$th urn dies and meanwhile gives birth to $k$ new particles in the $i$th urn at rate $\psi_k\left(\frac{i}{N}\right)$.

Let $X_t(i)$ be the number of particles in the $i$th urn at moment $t$ and
\[
X_t=\left(X_t(1), X_t(2),\ldots, X_t(N)\right),
\]
then our $N$-urn branching process $\{X_t\}_{t\geq 0}$ is a continuous-time Markov process with sate space
\[
\mathbb{X}=\left\{0,1,2,\ldots,\right\}^N.
\]
For any $x\in \mathbb{X}$ and $1\leq i\leq N$, we use $x(i)$ to denote the $i$th coordinate of $x$. For any $x\in \mathbb{X}, k\geq 0$ and $1\leq i\neq j\leq N$. we use $x^{(i,j)}_k$ to denote the configuration in $\mathbb{X}$ such that
\[
x^{(i,j)}_k(l)=
\begin{cases}
x(l) & \text{~if~}l\neq i,j,\\
x(i)-1 & \text{~if~}l=i,\\
x(j)+k & \text{~if~}l=j
\end{cases}
\]
for all $1\leq l\leq N$. For any $x\in \mathbb{X}, 1\leq i\leq N$ and integer $k$, we use $x^i_k$ to denote the configuration in $\mathbb{X}$ such that
\[
x^i_k(l)=
\begin{cases}
x(l) & \text{~if~}l\neq i,\\
x(i)+k & \text{~if~}l=i
\end{cases}
\]
for all $1\leq l\leq N$. According to the definition of $\{X_t\}_{t\geq 0}$, the generator $\mathcal{L}$ of $\{X_t\}_{t\geq 0}$ is given by
\begin{align}\label{equ 1.1 generator}
\mathcal{L}f(x)=&\frac{1}{N}\sum_{k=0}^{+\infty}\sum_{i=1}^N\sum_{j\neq i}^Nx(i)\lambda_k\left(\frac{i}{N}, \frac{j}{N}\right)\left(f(x^{(i,j)}_k)-f(x)\right)\notag\\
&+\sum_{k=0}^{+\infty}\sum_{i=1}^Nx(i)\psi_k\left(\frac{i}{N}\right)\left(f(x^i_{k-1})-f(x)\right)
\end{align}
for any $x\in \mathbb{X}$ and $f$ from $\mathbb{X}$ to $\mathbb{R}$. From now on, to distinguish different $N$, we write $X_t$ and $\mathcal{L}$ as $X_t^N$ and $\mathcal{L}_N$ respectively.

In this paper we assume that
\[
\sup_{u,v\in \mathbb{T}}\sum_{k\geq 0}k^{2+\varepsilon_0}\lambda_k(u,v), \text{~} \sup_{u\in \mathbb{T}}\sum_{k\geq 0}k^{2+\varepsilon_0}\psi_k(u)<+\infty
\]
and $\sum_{k\geq 0}k^{2+\varepsilon_0}\lambda_k\in C^\infty(\mathbb{T}^2), \sum_{k\geq 0}k^{2+\varepsilon_0}\psi_k\in C^\infty(\mathbb{T})$ for some $\varepsilon_0>0$.

\textbf{Example 1} Ehrenfest models. If $\lambda_k\equiv 0$ and $\psi_k\equiv 0$ for all $k\neq 1$, then our model reduces to the $N$-urn generalized Ehrenfest model (see \cite{Cheng2020} and \cite{Xue2020}), where particles perform independent random walks between $N$ urns such that a particle in the $i$th urn jumps to the $j$th one at rate $\frac{1}{N}\lambda_1\left(\frac{i}{N}, \frac{j}{N}\right)$.

\textbf{Example 2} Branching random walks. If $\lambda_k\equiv 0$ for all $k\neq 1$, then our model reduces to the $N$-urn branching random walk, which also includes the Ehrenfest model as a special case. In the branching random walk, each particle jumps between $N$-urns and eventually dies in a urn and meanwhile gives birth to some new particles in the same urn.

In \cite{Xue2020}, limit theorems of empirical density fields of $N$-urn generalized Ehrenfest models are given, which we will extend to the branching process case in this paper. That is to say, we will give law of large numbers and central limit theorem of $\sum_{i=1}^NX_t^N(i)\delta_{\frac{i}{N}}(du)$, where $\delta_a(du)$ is the Dirac measure concentrated on $a$. Inspired by Section 11.4 of \cite{Ethier1986}, we also consider limit theorems of hitting times of $\{X_t^N\}_{t\geq 0}$ as applications of our main results. For mathematical details, see Section \ref{section two}.

In the special case where $\lambda_k\equiv 0$ for all $k\neq 1$ and $\psi_k\equiv q_k$ for some $q_k\in [0, +\infty)$ for all $k\geq 0$, the total number $\mathcal{N}_t^N=\sum_{i=1}^NX_t^N(i)$ of particles is an example of density-dependent Markov chains introduced in \cite{Kurtz1978} such that
\[
\mathcal{N}_t^N\rightarrow \mathcal{N}_t^N+k-1 \text{~at rate~} N\beta_k\left(\frac{\mathcal{N}_t^N}{N}\right)
\]
for all $k\geq 0$, where $\beta_k(u)=q_k u$. Then limit theorems of $\mathcal{N}_t^N$ are applications of main theorems given in \cite{Kurtz1978} and Chapter 11 of \cite{Ethier1986}. Since $\mathcal{N}_t^N$ can be considered as $\sum_{i=1}^NX_t^N(i)\delta_{\frac{i}{N}}(f)$ for $f\equiv 1$, we can revisit limit theorems of $\mathcal{N}_t^N$ according to our main results. For mathematical details, see Section \ref{section six}.

\section{Main results}\label{section two}
In this section we give our main results. First we introduce some notations and definitions for later use. For any $a,b\in \mathbb{R}$, we use $\mathbb{N}(a, b^2)$ to denote the normal distribution with mean $a$ and variance $b^2$. For any integer $m\geq 1$, we use $\|\cdot\|_\infty$ to denote the $l_\infty$-norm on $C(\mathbb{T}^m)$, i.e.,
\[
\|f\|_\infty=\sup_{u\in \mathbb{T}^m}|f(u)|
\]
for any $f\in C(\mathbb{T}^m)$. We use $\|\cdot\|_2$ to denote the $l_2$-norm on $L^2(\mathbb{T}^m)$ for any integer $m\geq 1$, i.e.,
\[
\|f\|_2=\int_{T^m}f^2(u)du
\]
for any $f\in L^2(\mathbb{T}^m)$.

For any integer $m\geq 1$, we use $\mathcal{S}_m$ to denote the dual of $C^\infty(\mathbb{T}^m)$ endowed with weak topology, i.e., $\pi_n\rightarrow \pi$ in $\mathcal{S}_m$ if and only if
\[
\lim_{n\rightarrow+\infty}\pi_n(f)=\pi(f)
\]
for all $f\in C^\infty(\mathbb{T}^m)$. For simplicity, we write $\mathcal{S}_1$ as $\mathcal{S}$. For integers $m,k\geq 1$ and any continuous linear operator $\mathcal{R}$ from $C^\infty(\mathbb{T}^m)$ to $C^\infty(\mathbb{T}^k)$, we use $\mathcal{R}^*$ to denote the adjoint operator of $\mathcal{R}$, i.e., $\mathcal{R}^*$ is a continuous linear operator from $\mathcal{S}_k$ to $\mathcal{S}_m$ such that
\[
(\mathcal{R}^*\pi)(f)=\pi(\mathcal{R}f)
\]
for any $\pi\in \mathcal{S}_k$ and $f\in C^\infty(\mathbb{T}^m)$. For any $t\geq 0$, we define
\[
\mu_t^N(du)=\frac{1}{N}\sum_{i=1}^NX_t^N(i)\delta_{\frac{i}{N}}(du),
\]
where $\delta_a(du)$ is the Dirac measure concentrated on $a$, i.e., $\mu_t^N$ is in the dual of $C(\mathbb{T})$ and hence a $\mathcal{S}$-valued random variable such that
\[
\mu_t^N(f)=\frac{1}{N}\sum_{i=1}^NX_t^N(i)f\left(\frac{i}{N}\right)
\]
for any $f\in C(\mathbb{T})$. We define $\mathcal{P}_1, \mathcal{P}_2$ as continuous linear operators from $C(\mathbb{T})$ to $C(\mathbb{T})$ such that
\[
(\mathcal{P}_1f)(u)=\int_{\mathbb{T}}\sum_{k=0}^{+\infty}\lambda_k(u,v)(kf(v)-f(u))dv+f(u)\sum_{k=0}^{+\infty}(k-1)\psi_k(u)
\]
and
\[
(\mathcal{P}_2f)(u)=f(u)\sum_{k=0}^{+\infty}(k-1)\psi_k(u)+\int_{\mathbb{T}}\sum_{k=0}^{+\infty}k\lambda_k(v,u)f(v)dv
-f(u)\int_{\mathbb{T}}\sum_{k=0}^{+\infty}\lambda_k(u,v)dv
\]
for any $u\in\mathbb{T}$ and $f\in C(\mathbb{T})$. Note that $\mathcal{P}_1, \mathcal{P}_2$ are continuous since
\begin{equation}\label{equ 2.2 Lipschitz condition}
\|\mathcal{P}_if\|_\infty\leq \left(\left\|\sum_{k=0}^{+\infty}(k+1)\lambda_k\right\|_\infty+\left\|\sum_{k=0}^{+\infty}(k+1)\psi_k\right\|_\infty\right)\|f\|_\infty
\end{equation}
for any $f\in C(\mathbb{T})$. By Equation \eqref{equ 2.2 Lipschitz condition}, for any $t\geq 0$ and $i=1,2$ it is reasonable to define $e^{t\mathcal{P}_i}$ as
\[
e^{t\mathcal{P}_i}=\sum_{k=0}^{+\infty}\frac{t^k\mathcal{P}_i^k}{k!}.
\]
Since $\sum_{k\geq 0}k^{2+\varepsilon_0}\lambda_k\in C^\infty(\mathbb{T}^2)$ and $\sum_{k=0}^{+\infty}k^{2+\varepsilon_0}\psi_k\in C^\infty(\mathbb{T})$, it is easy to check that $\mathcal{P}_i$ and $e^{t\mathcal{P}_i}$ are also continuous linear operators from $C^\infty(\mathbb{T})$ to $C^\infty(\mathbb{T})$ for $i=1,2$.

Throughout this paper we adopt the following initial assumption.

\textbf{Assumption} (A): $\{X_0^N(i)\}_{1\leq i\leq N}$ are independent and $X_0^N(i)$ follows the Poisson distribution with mean $\phi\left(\frac{i}{N}\right)$ for all $1\leq i\leq N$, where $\phi\in C^\infty(\mathbb{T})$.

Now we can give our first main result about the law of large numbers of $\{\mu_t^N\}_{t\geq 0}$.

\begin{theorem}\label{theorem 2.1 scaling limit}
Under Assumption (A),
\[
\lim_{N\rightarrow+\infty}\mu_t^N(f)=\int_\mathbb{T}\rho_t(u)f(u)du
\]
in $L^2$ for any $t\geq 0$ and $f\in C(\mathbb{T})$, where $\{\rho_t\}_{t\geq 0}$ is the unique solution to the ordinary differential equation
\begin{equation}\label{equ 2.1 functionValuedODE}
\begin{cases}
&\frac{d}{dt}\rho_t=\mathcal{P}_2\rho_t,\\
&\rho_0=\phi
\end{cases}
\end{equation}
on $C(\mathbb{T})$ endowed with the norm $\|\cdot\|_\infty$.
\end{theorem}

\begin{remark}\label{remark 2.3}
By Equation \eqref{equ 2.2 Lipschitz condition}, Equation \eqref{equ 2.1 functionValuedODE} satisfies the global Lipschitz's condition. Hence, the unique solution $\{\rho_t\}_{t\geq 0}$ to Equation \eqref{equ 2.1 functionValuedODE} is given by
\[
\rho_t=e^{t\mathcal{P}_2}\phi
\]
according to Theorem 19.1.2 of \cite{Lang}. Furthermore, let $\mu_t(du)=\rho_t(u)du$, then $\{\mu_t\}_{t\geq 0}$ is the unique element from $[0, +\infty)$ to the dual of $C(\mathbb{T})$ such that
\[
\mu_t(f)=\int_\mathbb{T}\phi(u)f(u)du+\int_0^t \mu_s(\mathcal{P}_1f)ds
\]
for any $t\geq 0$ and $f\in C(\mathbb{T})$ according to the exchange of integration order.

\end{remark}

By Theorem \ref{theorem 2.1 scaling limit}, the scaling limit of $X_t^N$, i.e., law of large numbers of the empirical density field $\mu_t^N(du)$, is driven by a $C(\mathbb{T})$-valued linear ordinary differential equation. It is natural to further discuss the fluctuation of $X_t^N$, i.e., the central limit theorem from above law of large numbers of $\mu_t^N$. Hence, for any $t\geq 0$ we define
\[
V_t^N(du)=\frac{1}{\sqrt{N}}\sum_{i=1}^N\left(X_t^N(i)-\mathbb{E}X_t^N(i)\right)\delta_{\frac{i}{N}}(du),
\]
where $\mathbb{E}$ is the expectation operator. For some technical reasons we consider $V_t^N$ as random elements in $\mathcal{S}$, i.e.,
\[
V_t^N(f)=\frac{1}{\sqrt{N}}\sum_{i=1}^N\left(X_t^N(i)-\mathbb{E}X_t^N(i)\right)f\left(\frac{i}{N}\right)
\]
for all $f\in C^\infty(\mathbb{T})$. Then for any $T>0$, $\{V_t^N\}_{0\leq t\leq T}$ is a random element in $\mathcal{D}\left([0, T], \mathcal{S}\right)$, where $\mathcal{D}\left([0, T], \mathcal{S}\right)$ is the set of c\`{a}dl\`{a}g functions from $[0, T]$ to $\mathcal{S}$ endowed with the Skorokhod topology.

To state our main result about the fluctuation of $X_t^N$, we need introduce some notations and definitions as preliminaries. For any $s\geq 0$, we use $\mathcal{A}_s$ to denote the continuous linear operator from $C^\infty(\mathbb{T})$ to $C^\infty(\mathbb{T})$ such that
\[
(\mathcal{A}_sf)(u)=\sqrt{\sum_{k=1}^{+\infty}(k-1)^2\psi_k(u)}\sqrt{\rho_s(u)}f(u)
\]
for any $u\in \mathbb{T}$ and $f\in C^\infty(\mathbb{T})$, where $\rho_t$ is the unique solution to Equation \eqref{equ 2.1 functionValuedODE}. Note that $\rho_t\in C^\infty(\mathbb{T})$ since $e^{t\mathcal{P}_2}$ is from $C^\infty(\mathbb{T})$ to $C^\infty(\mathbb{T})$ and $\rho_t=e^{t\mathcal{P}_2}\phi$ as we have shown in Remark \ref{remark 2.3}. For any $s\geq 0$ and $k\geq 0$, we use $\mathcal{U}_s^k$ to denote the continuous linear operator from $C^\infty(\mathbb{T})$ to $C^\infty(\mathbb{T}^2)$ such that
\[
(\mathcal{U}^k_sf)(u,v)=\sqrt{\lambda_k(u,v)}\sqrt{\rho_s(u)}\left(kf(v)-f(u)\right)
\]
for any $u,v\in \mathbb{T}$ and $f\in C^\infty(\mathbb{T})$. We use $\{\mathcal{B}_s\}_{s\geq 0}$ to denote the standard Brownian motion on $L^2(\mathbb{T})$ such that $\{\mathcal{B}_s(f)\}_{s\geq 0}$ is a real-valued Brownian motion with
\[
{\rm Cov}\left(\mathcal{B}_s(f), \mathcal{B}_s(f)\right)=s\|f\|_2^2=s\int_\mathbb{T}f^2(u)du
\]
for any $s\geq 0$ and $f\in C^\infty(\mathbb{T})$. Similarly,  we use $\{\mathcal{W}_s\}_{s\geq 0}$ to denote the standard Brownian motion on $L^2(\mathbb{T}^2)$ such that $\{\mathcal{W}_s(f)\}_{s\geq 0}$ is a real-valued Brownian motion with
\[
{\rm Cov}\left(\mathcal{W}_s(f), \mathcal{W}_s(f)\right)=s\|f\|_2^2=s\int_{\mathbb{T}^2}f^2(u,v)dudv
\]
for any $s\geq 0$ and $f\in C^\infty(\mathbb{T}^2)$. Let $\{\mathcal{W}_t^k:~t\geq 0\}_{k\geq 0}$ be independent copies of $\{\mathcal{W}_t:~t\geq 0\}$ and be independent of $\{\mathcal{B}_t:~t\geq 0\}$, then for any $T>0$ we introduce the definition of the following $S$-valued Ornstein-Uhlenbeck (O-U) process $\{V_t\}_{0\leq t\leq T}$:
\begin{equation}\label{equ 2.3 GeneralizedOUprocess}
\begin{cases}
&dV_t=\mathcal{P}_1^*V_tdt+\mathcal{A}_t^*d\mathcal{B}_t+\sum_{k=0}^{+\infty}\left(\mathcal{U}_t^k\right)^*d\mathcal{W}_t^k \text{~for~}0\leq t\leq T,\\
&V_0\text{~is independent of~}\{\mathcal{B}_t\}_{t\geq 0}\text{~and~} \{\mathcal{W}_t^k:~t\geq 0\}_{k\geq 0}.
\end{cases}
\end{equation}
We define the solution to Equation \eqref{equ 2.3 GeneralizedOUprocess} as the solution to the martingale problem with respect to Equation \eqref{equ 2.3 GeneralizedOUprocess}. In detail, according to an analysis similar with that leading to Theorem 1.4 of \cite{Holley1978}, when the distribution of $V_0$ is given, there exists a unique random element $\{V_t\}_{0\leq t\leq T}$ in $\mathcal{D}\left([0, T], \mathcal{S}\right)$ with following two properties.

1) For any $f\in C^\infty(\mathbb{T})$, $\{V_t(f)\}_{0\leq t\leq T}$ is continuous in $t$.

2) For any $G\in C_c^\infty(\mathbb{R})$ and $f\in C^\infty(\mathbb{T})$,
\begin{align*}
\Bigg\{ G\left(V_t(f)\right)-G\left(V_0(f)\right)&-\int_0^tG^\prime\left(V_s(f)\right)V_s\left(\mathcal{P}_1f\right)ds\\
&-\frac{1}{2}\int_0^tG^{\prime\prime}\left(V_s(f)\right)\left(\|\mathcal{A}_sf\|_2^2+\sum_{k=0}^{+\infty}\|\mathcal{U}_s^kf\|_2^2\right)ds\Bigg\}_{0\leq t\leq T}
\end{align*}
is a martingale.

We define the above $\{V_t\}_{0\leq t\leq T}$ as the unique solution to Equation \eqref{equ 2.3 GeneralizedOUprocess}. Now we can give our second main result, which is about the fluctuation of $X_t^N$.

\begin{theorem}\label{theorem 2.2 fluctuation}
Under Assumption (A), $\{V_t^N\}_{0\leq t\leq T}$ converges weakly to $\{V_t\}_{0\leq t\leq T}$ as $N\rightarrow+\infty$ for any $T>0$, where $\{V_t\}_{0\leq t\leq T}$ is the unique solution to Equation \eqref{equ 2.3 GeneralizedOUprocess} with $V_0(f)$ following normal distribution $\mathbb{N}\left(0, \int_{\mathbb{T}}\phi(u)f^2(u)du\right)$ for all $f\in C^\infty(\mathbb{T})$.
\end{theorem}

\begin{remark}\label{remark 2.4}
Since $e^{t\mathcal{P}_1}$ is a continuous linear operator from $C^\infty(\mathbb{T})$ to $C^\infty(\mathbb{T})$, $(e^{t\mathcal{P}_1})^*$ is well-defined. Furthermore,
\[
(e^{t\mathcal{P}_1})^*=\sum_{k=0}^{+\infty}\frac{t^k(\mathcal{P}_1^*)^k}{k!},
\]
which can be formally written as $e^{t\mathcal{P}_1^*}$. Hence it is also reasonable to define the solution $V_t$ to Equation \eqref{equ 2.3 GeneralizedOUprocess} as
\begin{equation}\label{equ 2.4}
V_t=e^{t\mathcal{P}_1^*}V_0+\int_0^te^{(t-s)\mathcal{P}_1^*}\mathcal{A}_s^*d\mathcal{B}_t
+\sum_{k=0}^{+\infty}\int_0^te^{(t-s)\mathcal{P}_1^*}(\mathcal{U}_t^k)^*d\mathcal{W}_t^k.
\end{equation}
An important and natural question is whether above two definitions of $V_t$ are equivalent. The answer is positive since the analysis leading to the uniqueness of the solution to the martingale problem gives Fourier transforms of the distribution of this solution, which coincide with Fourier transforms of the distribution of $V_t$ given by Equation \eqref{equ 2.4}. For mathematical details, see the proof of Theorem 1.4 of \cite{Holley1978} or a similar analysis given in Appendix A.2 of \cite{Xue2020}.
\end{remark}

By Theorem \ref{theorem 2.2 fluctuation}, the fluctuation of $X_t^N$ is driven by a $\mathcal{S}$-valued O-U process. By Remark \ref{remark 2.4}, we have the following corollary of Theorem \ref{theorem 2.2 fluctuation}.

\begin{corollary}\label{corollary 2.3}
Under Assumption (A), $V_t^N(f)$ converges in distribution to $\mathbb{N}\left(0, \theta_t^2(f)\right)$ as $N\rightarrow+\infty$ for any $t\geq 0$ and $f\in C^\infty(\mathbb{T})$, where
\[
\theta_t^2(f)=\int_{\mathbb{T}}\phi(u)\left(e^{t\mathcal{P}_1}f\right)^2(u)du
+\int_0^t \left\|\mathcal{A}_se^{(t-s)\mathcal{P}_1}f\right\|_2^2ds+\sum_{k=0}^{+\infty}\int_0^t\left\|\mathcal{U}_s^ke^{(t-s)\mathcal{P}_1}f\right\|_2^2ds.
\]
\end{corollary}

\begin{remark}\label{remark new difficulty}
Theorems \ref{theorem 2.1 scaling limit} and \ref{theorem 2.2 fluctuation} extend results about scaling limit and fluctuation of the $N$-urn Ehrenfest model to the branching process case. The main new difficulty in the branching process case is that the total number of particles is not conserved, which makes several techniques given in \cite{Xue2020} for the Ehrenfest model not feasible for the branching process case. For example, it is easy to show that numbers of particles in different urns are uncorrelated in the Ehrenfest model by introducing fixed number of indicator functions recording positions of every particles. However, in the branching process case new particles are born randomly and then the proof of approximated independence between numbers of particles in different urns follows from a totally different approach, where covariance of numbers of particles in two different urns is bounded by the difference of respective solutions to two different linear ordinary differential equations. For mathematical details, see Section \ref{section three}.

\end{remark}

As applications of our main result, we discuss limit theorems of hitting times of $\mu_t^N$. For given $f\in C^\infty(\mathbb{T})$ and $r>\int_\mathbb{T}f(u)\phi(u)du$, we define
\[
\tau^N_{r,f}=\inf\left\{t:~\mu_t^N(f)\geq r\right\}.
\]
Let $\mu_t(du)=\rho_t(u)du$ as we have defined in Remark \ref{remark 2.3}, then we define
\[
\tau_{r,f}=\inf\left\{t:~\mu_t(f)=r\right\}.
\]
We have the following result about law of large numbers and central limit theorem of $\tau^N_{r,f}$.

\begin{theorem}\label{theorem 2.3 hittingTimes}
If $\tau_{r,f}<+\infty$ and $\mu_{\tau_{r,f}}\left(\mathcal{P}_1f\right)>0$, then $\tau_{r,f}^N$ converges to $\tau_{r,f}$ in probability and
\[
\sqrt{N}\left(\tau^N_{r,f}-\tau_{r,f}\right)
\]
converges to $\mathbb{N}\left(0, \frac{\theta^2_{\tau_{r,f}}(f)}{\mu^2_{\tau_{r,f}}(\mathcal{P}_1f)}\right)$ in distribution as $N\rightarrow+\infty$.
\end{theorem}

\begin{remark}\label{remark 2.5}
As we have introduced in Remark \ref{remark 2.3},
\[
\frac{d}{dt}\mu_t(f)=\mu_t(\mathcal{P}_1f).
\]
Hence, if $f\in C^\infty(\mathbb{T})$ makes $\mu_0(\mathcal{P}_1f)>0$, then there exists $T_0>0$ such that $\mu_s(\mathcal{P}_1f)>0$ for $s\in [0, T_0]$ and $\mu_t(f)$ is increasing in $t\in [0, T_0]$. Let $r\in (\mu_0(f), \mu_{T_0}(f))$, then $(f, r)$ satisfies the assumption in Theorem \ref{theorem 2.3 hittingTimes}.
\end{remark}

The remainder of this paper is organized as follows. As a preliminary for proofs of main theorems, in Section \ref{section three} we show that numbers of particles in different urns are approximately independent as $N\rightarrow+\infty$, which is the main difficulty in our proofs. Proofs of Theorems \ref{theorem 2.1 scaling limit} and \ref{theorem 2.2 fluctuation} are given in Section \ref{section four}. With the approximate independence between numbers of particles in different urns given in Section \ref{section three}, Theorem \ref{theorem 2.1 scaling limit} follows from a mean-variance analysis and Theorem \ref{theorem 2.2 fluctuation} follows from a martingale strategy. Proof of Theorem \ref{theorem 2.3 hittingTimes} is given in Section \ref{section five}. According to the fact that $\mu^N_{\tau^N_{r,f}}(f)\approx \mu_{\tau_{r,f}}(f)$, the central limit of $\tau_{r,f}^N$ can be related to fluctuation of $\{\mu_t^N\}_{t\geq 0}$. As applications of our main results, in Section \ref{section six} we revisit a special case which reduces to a density-dependent Markov chain. We apply our main results to give new proofs of limit theorems of total particles number which are corollaries of main theorems given in \cite{Kurtz1978} and Chapter 11 of \cite{Ethier1986}.

\section{Approximate independence}\label{section three}
In this section we prove the following lemma.
\begin{lemma}\label{lemma 3.1 approximate independence}
For any $t\geq 0$, there exists $C_1=C_1(t)<+\infty$ independent of $N$ such that
\[
\left|{\rm Cov}\left(X_s^N(i), X_s^N(j)\right)\right|\leq \frac{C_1}{N}
\]
for all $N\geq 1$, $1\leq i\neq j\leq N$ and $0\leq s\leq t$.
\end{lemma}

For the Ehrenfest model case, it is not difficult to show that ${\rm Cov}\left(X_s^N(i), X_s^N(j)\right)=0$ for $i\neq j$ under Assumption (A) according to the fact that $\sum_{i=0}^{+\infty}X_t^N(i)=\sum_{i=0}^{+\infty}X_0^N(i)$. For the branching process case, new particles are born randomly and hence we only have Lemma \ref{lemma 3.1 approximate independence}, which is a weaker conclusion that covariance of numbers of particles in two different urns is $O(N^{-1})$.

As a preliminary for the proof of Lemma \ref{lemma 3.1 approximate independence}, we introduce some definitions and notations. For a finite set $A$, we call a function from $A^2$ to $\mathbb{R}$ a $A\times A$ matrix. For any two $A\times A$ matrices $M_1, M_2$, we define $M_1M_2$ as the $A\times A$ matrix such that
\[
M_1M_2(k,l)=\sum_{m\in A}M_1(k,m)M_2(m,l)
\]
for any $k,l\in A$. Then for a $A\times A$ matrix $M$, $t\geq 0$ and any inter $n\geq 2$, we define $M^2=MM$, $M^n=M^{n-1}M$ for $n\geq 3$ and
\[
e^{tM}=\sum_{k=0}^{+\infty}\frac{t^kM^k}{k!}.
\]
Furthermore, for a $A\times A$ matrix $M$ and function $H$ from $A$ to $\mathbb{R}$, we define $MH$ as the function from $A$ to $\mathbb{R}$ such that
\[
MH(k)=\sum_{m\in A}M(k,m)H(m)
\]
for all $k\in A$. Let $\mathbb{T}^N=\{1,2,\ldots, N\}$, then we define $\Gamma_t^N, \hat{\Gamma}_t^N$ as functions from $(\mathbb{T}^N)^2$ to $\mathbb{R}$ such that
\[
\Gamma_t^N(i,j)=\mathbb{E}X_t^N(i)\mathbb{E}X_t^N(j)\text{~and~}\hat{\Gamma}_t^N(i,j)=\mathbb{E}\left(X_t^N(i)X_t^N(j)\right)
\]
for all $i,j\in \mathbb{T}^N$. As a result,
\[
{\rm Cov}\left(X_t^N(i), X_t^N(j)\right)=\hat{\Gamma}_t^N(i,j)-\Gamma_t(i,j)
\]
for all $i,j\in \mathbb{T}^N$. Now we give the proof of Lemma \ref{lemma 3.1 approximate independence}.

\proof[Proof of Lemma \ref{lemma 3.1 approximate independence}]

According to the generator $\mathcal{L}_N$ of $\{X_t^N\}_{t\geq 0}$ given in Section \ref{section one} and Chapman-Kolmogorov equation, for $t\geq 0$ and $1\leq i\leq N$,
\begin{align}\label{equ 3.1}
\frac{d}{dt}\mathbb{E}X_t^N(i)=&-\frac{\mathbb{E}X_t^N(i)}{N}\sum_{k=0}^{+\infty}\sum_{j\neq i}\lambda_k\left(\frac{i}{N}, \frac{j}{N}\right)
+\frac{1}{N}\sum_{k=0}^{+\infty}\sum_{j\neq i}k\lambda_k\left(\frac{j}{N},\frac{i}{N}\right)\mathbb{E}X_t^N(j) \notag\\
&+\mathbb{E}X_t^N(i)\sum_{k=0}^{+\infty}(k-1)\psi_k\left(\frac{i}{N}\right).
\end{align}
As a result,
\[
\frac{d}{dt}\Gamma_t^N=\mathcal{M}^N_1\Gamma_t^N
\]
and hence $\Gamma_t^N=e^{t\mathcal{M}^N_1}\Gamma_0^N$, where $\mathcal{M}^N_1$ is a $(\mathbb{T}^N)^2\times (\mathbb{T}^N)^2$ matrix such that
\begin{align*}
&\mathcal{M}^N_1\left((i,j), (w,m)\right)=\\
&\begin{cases}
&\frac{1}{N}\sum_{k=0}^{+\infty}k\lambda_k\left(\frac{w}{N}, \frac{i}{N}\right) \text{\quad if~}i\neq j, w\neq i\text{~and~}m=j,\\
&\frac{1}{N}\sum_{k=0}^{+\infty}k\lambda_k\left(\frac{m}{N}, \frac{j}{N}\right) \text{\quad if~}i\neq j, w=i\text{~and~}m\neq j,\\
&-\frac{1}{N}\sum_{k=0}^{+\infty}\sum_{l\neq i}\lambda_k\left(\frac{i}{N}, \frac{l}{N}\right)
-\frac{1}{N}\sum_{k=0}^{+\infty}\sum_{l\neq j}\lambda_k\left(\frac{j}{N}, \frac{l}{N}\right)\\
&\text{\quad}+\sum_{k=0}^{+\infty}(k-1)\psi_k\left(\frac{i}{N}\right)+\sum_{k=0}^{+\infty}(k-1)\psi_k\left(\frac{j}{N}\right)
\text{\quad if~}i\neq j\text{~and~}(w,m)=(i,j),\\
&-\frac{2}{N}\sum_{k=0}^{+\infty}\sum_{l\neq i}\lambda_k\left(\frac{i}{N}, \frac{l}{N}\right)
+2\sum_{k=0}^{+\infty}(k-1)\psi_k\left(\frac{i}{N}\right) \text{\quad if ~}i=j \text{~and~} (w,m)=(i,i),\\
&\frac{1}{N}\sum_{k=0}^{+\infty}k\lambda_k\left(\frac{w}{N}, \frac{i}{N}\right) \text{\quad if~}i=j, w\neq i\text{~and~}m=i,\\
&\frac{1}{N}\sum_{k=0}^{+\infty}k\lambda_k\left(\frac{m}{N}, \frac{i}{N}\right) \text{\quad if~}i=j, w=i\text{~and~}m\neq i,\\
&0\text{\quad else}.
\end{cases}
\end{align*}
Similarly, according to the generator $\mathcal{L}_N$ and Chapman-Kolmogorov equation,
\begin{align*}
&\frac{d}{dt}\mathbb{E}\left(X_t^N(i)X_t^N(j)\right)\\
&=\sum_{k=0}^{+\infty}\frac{\lambda_k\left(\frac{i}{N}, \frac{j}{N}\right)}{N}
\mathbb{E}\left(X_t^N(i)\left(-X_t^N(j)+kX_t^N(i)-k\right)\right)\\
&\text{\quad}+\sum_{k=0}^{+\infty}\frac{\lambda_k\left(\frac{j}{N},\frac{i}{N}\right)}{N}\mathbb{E}\left(X_t^N(j)\left(-X_t^N(i)+kX_t^N(j)-k\right)\right)\\
&\text{\quad}-\frac{\mathbb{E}\left(X_t^N(i)X_t^N(j)\right)}{N}\sum_{k=0}^{+\infty}\sum_{l\neq i,j}\lambda_k\left(\frac{i}{N}, \frac{l}{N}\right)
-\frac{\mathbb{E}\left(X_t^N(i)X_t^N(j)\right)}{N}\sum_{k=0}^{+\infty}\sum_{l\neq i,j}\lambda_k\left(\frac{j}{N}, \frac{l}{N}\right)\\
&\text{\quad}+\sum_{k=0}^{+\infty}\sum_{l\neq i,j}\frac{k\mathbb{E}\left(X_t^N(l)X_t^N(j)\right)}{N}\lambda_k\left(\frac{l}{N}, \frac{i}{N}\right)
+\sum_{k=0}^{+\infty}\sum_{l\neq i,j}\frac{k\mathbb{E}\left(X_t^N(i)X_t^N(l)\right)}{N}\lambda_k\left(\frac{l}{N}, \frac{j}{N}\right)\\
&\text{\quad}+\sum_{k=0}^{+\infty}(k-1)\mathbb{E}\left(X_t^N(i)X_t^N(j)\right)\psi_k\left(\frac{i}{N}\right)
+\sum_{k=0}^{+\infty}(k-1)\mathbb{E}\left(X_t^N(i)X_t^N(j)\right)\psi_k\left(\frac{j}{N}\right)
\end{align*}
for $i\neq j$ and
\begin{align*}
&\frac{d}{dt}\mathbb{E}\left(\left(X_t^N(i)\right)^2\right)\\
&=\frac{1}{N}\sum_{l\neq i}\sum_{k=0}^{+\infty}\lambda_k\left(\frac{i}{N}, \frac{l}{N}\right)\mathbb{E}\left(X_t^N(i)\left(-2X_t^N(i)+1\right)\right)\\
&\text{\quad}+\frac{1}{N}\sum_{l\neq i}\sum_{k=0}^{+\infty}\lambda_k\left(\frac{l}{N}, \frac{i}{N}\right)
\mathbb{E}\left(X_t^N(l)\left(2kX_t^N(i)+k^2\right)\right)\\
&\text{\quad}+\sum_{k=0}^{+\infty}\psi_k\left(\frac{i}{N}\right)\mathbb{E}\left(X_t^N(i)\left(2(k-1)X_t^N(i)+(k-1)^2\right)\right).
\end{align*}
Then, since $X_t^N(i)\leq \left(X_t^N(i)\right)^2$ for all $1\leq i\leq N$,
\[
\mathcal{M}_2^N\hat{\Gamma}_t^N(i,j)\leq \frac{d}{dt}\hat{\Gamma}_t^N(i,j)\leq \mathcal{M}_3^N\hat{\Gamma}_t^N(i,j)
\]
for any $1\leq i,j\leq N$, where $\mathcal{M}_2^N, \mathcal{M}_3^N$ are $(\mathbb{T}^N)^2\times (\mathbb{T}^N)^2$ matrices such that
\begin{align*}
&\mathcal{M}^N_2\left((i,j), (w,m)\right)=\\
&\begin{cases}
&\frac{1}{N}\sum_{k=0}^{+\infty}k\lambda_k\left(\frac{w}{N}, \frac{i}{N}\right) \text{\quad if~}i\neq j, w\neq i,j\text{~and~}m=j,\\
&\frac{1}{N}\sum_{k=0}^{+\infty}k\lambda_k\left(\frac{m}{N}, \frac{j}{N}\right) \text{\quad if~}i\neq j, w=i\text{~and~}m\neq i,j,\\
&-\frac{1}{N}\sum_{k=0}^{+\infty}\sum_{l\neq i}\lambda_k\left(\frac{i}{N}, \frac{l}{N}\right)
-\frac{1}{N}\sum_{k=0}^{+\infty}\sum_{l\neq j}\lambda_k\left(\frac{j}{N}, \frac{l}{N}\right)\\
&\text{\quad}+\sum_{k=0}^{+\infty}(k-1)\psi_k\left(\frac{i}{N}\right)+\sum_{k=0}^{+\infty}(k-1)\psi_k\left(\frac{j}{N}\right)
\text{\quad if~}i\neq j\text{~and~}(w,m)=(i,j),\\
&-\frac{2}{N}\sum_{k=0}^{+\infty}\sum_{l\neq i}\lambda_k\left(\frac{i}{N}, \frac{l}{N}\right)
+2\sum_{k=0}^{+\infty}(k-1)\psi_k\left(\frac{i}{N}\right) \text{\quad if ~}i=j \text{~and~} (w,m)=(i,i),\\
&\frac{1}{N}\sum_{k=0}^{+\infty}k\lambda_k\left(\frac{w}{N}, \frac{i}{N}\right) \text{\quad if~}i=j, w\neq i\text{~and~}m=i,\\
&\frac{1}{N}\sum_{k=0}^{+\infty}k\lambda_k\left(\frac{m}{N}, \frac{i}{N}\right) \text{\quad if~}i=j, w=i\text{~and~}m\neq i,\\
&0\text{\quad else}.
\end{cases}
\end{align*}
and
\begin{align*}
&\mathcal{M}^N_3\left((i,j), (w,m)\right)=\\
&\begin{cases}
&\frac{1}{N}\sum_{k=0}^{+\infty}k\lambda_k\left(\frac{w}{N}, \frac{i}{N}\right) \text{\quad if~}i\neq j, w\neq i\text{~and~}m=j,\\
&\frac{1}{N}\sum_{k=0}^{+\infty}k\lambda_k\left(\frac{m}{N}, \frac{j}{N}\right) \text{\quad if~}i\neq j, w=i\text{~and~}m\neq j,\\
&-\frac{1}{N}\sum_{k=0}^{+\infty}\sum_{l\neq i}\lambda_k\left(\frac{i}{N}, \frac{l}{N}\right)
-\frac{1}{N}\sum_{k=0}^{+\infty}\sum_{l\neq j}\lambda_k\left(\frac{j}{N}, \frac{l}{N}\right)\\
&\text{\quad}+\sum_{k=0}^{+\infty}(k-1)\psi_k\left(\frac{i}{N}\right)+\sum_{k=0}^{+\infty}(k-1)\psi_k\left(\frac{j}{N}\right)
\text{\quad if~}i\neq j\text{~and~}(w,m)=(i,j),\\
&-\frac{1}{N}\sum_{k=0}^{+\infty}\sum_{l\neq i}\lambda_k\left(\frac{i}{N}, \frac{l}{N}\right)\\
&\text{\quad}+\sum_{k=0}^{+\infty}\left(2(k-1)+(k-1)^2\right)\psi_k\left(\frac{i}{N}\right) \text{\quad if ~}i=j \text{~and~} (w,m)=(i,i),\\
&\frac{1}{N}\sum_{k=0}^{+\infty}k\lambda_k\left(\frac{w}{N}, \frac{i}{N}\right) \text{\quad if~}i=j, w\neq i\text{~and~}m=i,\\
&\frac{1}{N}\sum_{k=0}^{+\infty}k\lambda_k\left(\frac{m}{N}, \frac{i}{N}\right) \text{\quad if~}i=j, w=i\text{~and~}m\neq i,\\
&\frac{1}{N}\sum_{k=0}^{+\infty}k^2\lambda_k\left(\frac{w}{N}, \frac{i}{N}\right)\text{\quad if~}i=j \text{~and~} w=m\neq i,\\
&0\text{\quad else}.
\end{cases}
\end{align*}
Since $\mathcal{M}_l^N\left((i,j), (w,m)\right)\geq 0$ when $(w,m)\neq (i,j)$ for $l=2,3$, we have
\begin{equation}\label{equ 3.2}
\underline{\Lambda}_t^N(i,j)\leq \hat{\Gamma}_t^N(i,j)\leq \overline{\Lambda}_t^N(i,j)
\end{equation}
for all $i,j\in \{1,2,\ldots, N\}$, where $\overline{\Lambda}_t^N$ is the solution to
\[
\begin{cases}
&\frac{d}{dt}\overline{\Lambda}_t^N=\mathcal{M}_3^N\overline{\Lambda}_t^N,\\
&\overline{\Lambda}_0^N=\hat{\Gamma}_0^N
\end{cases}
\]
and $\underline{\Lambda}_t^N$ is the solution to
\[
\begin{cases}
&\frac{d}{dt}\underline{\Lambda}_t^N=\mathcal{M}_2^N\underline{\Lambda}_t^N,\\
&\underline{\Lambda}_0^N=\hat{\Gamma}_0^N,
\end{cases}
\]
i.e., $\overline{\Lambda}_t^N=e^{t\mathcal{M}_3^N}\hat{\Gamma}_0^N$ and $\underline{\Lambda}_t^N=e^{t\mathcal{M}_2^N}\hat{\Gamma}_0^N$. As a result, for $i, j\in \{1,2,\ldots, N\}$ and any $t\geq 0$,
\begin{align}\label{equ 3.3 compare}
&-\left|e^{t\mathcal{M}_2^N}\hat{\Gamma}_0^N(i,j)-e^{t\mathcal{M}_1^N}\Gamma_0^N(i,j)\right| \\
&\leq {\rm Cov}\left(X_t^N(i), X_t^N(j)\right)\leq \left|e^{t\mathcal{M}_3^N}\hat{\Gamma}_0^N(i,j)-e^{t\mathcal{M}_1^N}\Gamma_0^N(i,j)\right|. \notag
\end{align}
To bound lower and upper bounds in the above inequality, we point out some properties of $\mathcal{M}_l^N$ for $l=1,2,3$.
According to definitions of $\mathcal{M}_l^N$ for $l=1,2,3$, there exists $C_2<+\infty$ independent of $N$ such that
\begin{equation}\label{equ 3.4}
|\mathcal{M}_l^NH(i,j)|\leq C_2|H(i,j)|
\end{equation}
for any $H$ from $(\mathbb{T}^N)^2$ to $\mathbb{R}$, $l=1,2,3$ and $1\leq i,j\leq N$. For any $i\neq j$,
\begin{equation}\label{equ 3.5}
{\rm card}\left\{(w,m):~\mathcal{M}_l^N\left((i,j), (w,m)\right)\neq 0\right\}\leq 3N
\end{equation}
for $l=1,2,3$, where ${\rm card}(A)$ is the cardinality of set $A$. For any $i\neq j$ and $l=2,3$,
\begin{equation}\label{equ 3.6}
{\rm card}\left\{(w,m):~\mathcal{M}_l^N\left((i,j), (w,m)\right)\neq \mathcal{M}_1^N\left((i,j), (w,m) \right)\right\}\leq 2.
\end{equation}
For any $i\neq j$,
\begin{equation}\label{equ 3.7}
\mathcal{M}_1^N\left((i,j), (i,j)\right)=\mathcal{M}_2^N\left((i,j), (i,j)\right)=\mathcal{M}_3^N\left((i,j), (i,j)\right).
\end{equation}
There exists $C_3<+\infty$ independent of $N$ such that
\begin{equation}\label{equ 3.8}
\left|\mathcal{M}_l^N\left((i,j), (w,m)\right)\right|\leq \frac{C_3}{N}
\end{equation}
for any $1\leq i,j\leq N, (w,m)\neq (i,j)$ and $l=1,2,3$. There exists $C_4<+\infty$ independent of $N$ such that
\begin{equation}\label{equ 3.9}
\left|\mathcal{M}_l^N\left((i,j), (i,j)\right)\right|\leq C_4
\end{equation}
for any $1\leq i,j\leq N$.

Now we bound the upper bound in inequality \eqref{equ 3.3 compare} from above. For integer $k\geq 0$, let
\[
\zeta_k^N=\sup\left\{\left|\left(\mathcal{M}_3^N\right)^k\hat{\Gamma}^N_0(i,j)
-\left(\mathcal{M}_1^N\right)^k\Gamma^N_0(i,j)\right|:~1\leq i\neq j\leq N\right\},
\]
then $\zeta_0^N=0$ according to Assumption (A). Since
\[
\left(\mathcal{M}_3^N\right)^{k+1}\hat{\Gamma}^N_0(i,j)=\sum_{(w,m)}\mathcal{M}_3^N\left((i,j), (w,m)\right)\left(\mathcal{M}_3^N\right)^{k}\hat{\Gamma}^N_0(w,m)
\]
and
\[
\left(\mathcal{M}_1^N\right)^{k+1}\Gamma^N_0(i,j)=\sum_{(w,m)}\mathcal{M}_1^N\left((i,j), (w,m)\right)\left(\mathcal{M}_1^N\right)^{k}\Gamma^N_0(w,m),
\]
by Equations \eqref{equ 3.4} to \eqref{equ 3.9} we have
\begin{align*}
\zeta_{k+1}^N\leq C_4\zeta_k^N+3N\frac{C_3}{N}\zeta_k^N+\frac{4C_3}{N}C_2^k\left(\|\phi\|_\infty^2+\|\phi\|_\infty\right)
=C_5 \zeta_k^N+\frac{C_6C_2^k}{N},
\end{align*}
where $C_5=C_4+3C_3$ and $C_6=4C_3\left(\|\phi\|_\infty^2+\|\phi\|_\infty\right)$. Since $\zeta_0^N=0$, by induction,
\[
\zeta_k^N\leq \frac{C_6\sum_{m=0}^{k-1}C_2^mC_5^{k-1-m}}{N}\leq \frac{C_6(C_2+C_5)^{k-1}}{N}
\]
for $k\geq 1$. As a result, for $i\neq j$,
\[
\left|e^{t\mathcal{M}_3^N}\hat{\Gamma}_0^N(i,j)-e^{t\mathcal{M}_1^N}\Gamma_0^N(i,j)\right|\leq \frac{C_6t}{N}e^{t(C_2+C_5)}.
\]
According to a similar analysis, the above inequality still holds when we replace $\mathcal{M}_3^N$ by $\mathcal{M}_2^N$. Therefore, Lemma \ref{lemma 3.1 approximate independence} holds with $C_1=C_1(t)=C_6te^{t(C_2+C_5)}$.

\qed

\section{Proofs of Theorems \ref{theorem 2.1 scaling limit} and \ref{theorem 2.2 fluctuation}}\label{section four}
In this section we prove Theorems \ref{theorem 2.1 scaling limit} and \ref{theorem 2.2 fluctuation}.

\proof[Proof of Theorem \ref{theorem 2.1 scaling limit}]

By Equation \eqref{equ 3.4},
\begin{equation}\label{equ 4.1}
\mathbb{E}X_t^N(i)=\left(\Gamma_t^N(i,i)\right)^{\frac{1}{2}}\leq e^{\frac{1}{2}C_2t}\|\phi\|_\infty
\end{equation}
and
\begin{align}\label{equ 4.2}
{\rm Var}\left(X_t^N(i)\right)&\leq \hat{\Gamma}_t^N(i,i)\leq \sup_{1\leq i,j\leq N}\hat{\Gamma}_t^N(i,j)\\
&=\sup_{1\leq i,j\leq N}\mathbb{E}\left(X_t^N(i)X_t^N(j)\right)\leq e^{C_2t}\left(\|\phi\|_\infty^2+\|\phi\|_\infty\right) \notag
\end{align}
for any $1\leq i\leq N$ and $t\geq 0$. Hence, for $f\in C(\mathbb{T})$,
\begin{align*}
&{\rm Var}\left(\mu_t^N(f)\right)\leq\frac{\|f\|^2_\infty}{N^2}\sum_{i=1}^N{\rm Var}\left(X_t^N(i)\right)+\frac{\|f\|_\infty^2}{N^2}\sum_{i\neq j}\left|{\rm Cov}\left(X_t^N(i), X_t^N(j)\right)\right|\\
&\leq \frac{e^{C_2t}\left(\|\phi\|_\infty^2+\|\phi\|_\infty\right)\|f\|_\infty^2}{N}+\frac{C_1(t)\|f\|_\infty^2}{N}\rightarrow 0
\end{align*}
as $N\rightarrow+\infty$ by Equation \eqref{equ 4.2} and Lemma \ref{lemma 3.1 approximate independence}. As a result, by Chebyshev's inequality, to prove Theorem \ref{theorem 2.1 scaling limit} we only need to show that
\begin{equation}\label{equ 4.3}
\lim_{N\rightarrow+\infty}\varpi_t^N(f)=\mu_t(f)
\end{equation}
for any $f\in C(\mathbb{T})$, where $\mu_t(du)=\rho_t(u)du$ as we have defined in Section \ref{section two} and
\[
\varpi_t^N(du)=\frac{1}{N}\sum_{i=1}^N\mathbb{E}X_t^N(i)\delta_{\frac{i}{N}}(du).
\]
For $u\in \mathbb{T}$ and $t\geq 0$, we define $\rho_t^N(u)=\mathbb{E}X_t^N(i)$ when $\frac{i-1}{N}<u\leq \frac{i}{N}$ for some $i\in \{1,2,\ldots, N\}$. By Equations \eqref{equ 3.1} and \eqref{equ 4.1},
\[
\sup_{0\leq s\leq t, u\in \mathbb{T}}\left|\frac{d}{ds}\rho_s^N(u)-\mathcal{P}_2\rho_s^N(u)\right|=O(N^{-1}).
\]
Since $\phi\in C^\infty(\mathbb{T})$, $\|\rho_0-\rho_0^N\|_\infty=O(N^{-1})$. As a result, according to the fact that $\frac{d}{dt}\rho_t(u)=\mathcal{P}_2\rho_t(u)$,
\[
\|\rho_s^N-\rho_s\|_\infty\leq O(N^{-1})+\int_0^s \|\mathcal{P}_2(\rho_r^N-\rho_r^N)\|_\infty dr
\]
for $0\leq s\leq t$, where $O(N^{-1})$ can be chosen uniformly in $s\in [0, t]$. Since
\[
\|\mathcal{P}_2f\|_\infty\leq \left(\left\|\sum_{k=0}^{+\infty}(k+1)\psi_k\right\|_\infty+\left\|\sum_{k=0}^{+\infty}(k+1)\lambda_k\right\|_\infty\right)\|f\|_\infty
\]
for any $f\in C(\mathbb{T})$, by Grownwall's inequality we have
\begin{equation}\label{equ 4.4}
\|\rho_t^N-\rho_t\|_\infty \leq O(N^{-1})\exp\left\{t\left(\left\|\sum_{k=0}^{+\infty}(k+1)\psi_k\right\|_\infty+\left\|\sum_{k=0}^{+\infty}(k+1)\lambda_k\right\|_\infty\right)\right\}.
\end{equation}
Equation \eqref{equ 4.3} follows from Equation \eqref{equ 4.4} and the proof is complete.

\qed

As a preliminary for the proof of Theorem \ref{theorem 2.2 fluctuation}, we need following two lemmas.

\begin{lemma}\label{lemma 4.1}
For any $T\geq 0$, $\{V_t^N:~0\leq t\leq T\}_{N\geq 1}$ are tight.
\end{lemma}

\begin{lemma}\label{lemma 4.2}
For any $\epsilon>0$ and $t>0$,
\[
\lim_{N\rightarrow+\infty}P\left(\sup_{0\leq s\leq t}\left|\mu_s^N(f)-\mu_{s-}^N(f)\right|\geq \frac{\epsilon}{N^{\frac{4+3\varepsilon_0}{8+4\varepsilon_0}}}\right)=0,
\]
where $s-$ is the moment just before $s$, i.e., $\mu_{s-}^N(f)=\lim_{r<s,r\rightarrow s}\mu_r^N(f)$.
\end{lemma}

Lemma \ref{lemma 4.1} ensures that $\{V_t^N:~0\leq t\leq T\}_{N\geq 1}$ has weakly convergent subsequence. Lemma \ref{lemma 4.2} ensures that any weak limits of $\{\mu_s^N(f):~0\leq s\leq t\}_{N\geq 1}$ and $\{V_s^N(f):~0\leq s\leq t\}_{N\geq 1}$ are continuous.
This lemma is trivial in the Ehrenfest case where $\mu_s^N(f)-\mu_{s-}^N(f)=O(N^{-1})$.

\proof[Proof of Lemma \ref{lemma 4.1}]

By Aldous' criterion and Theorem 4.1 of \cite{Mitoma1983}, Lemma \ref{lemma 4.1} follows from following two Equations.

(1) For any $t>0$ and $f\in C^\infty(\mathbb{T})$,
\begin{equation}\label{equ 4.8}
\lim_{M\rightarrow+\infty}\limsup_{N\rightarrow+\infty}P\big(|V_t^N(f)|\geq M\big)=0.
\end{equation}

(2) For any $\epsilon>0$ and $f\in C^\infty(\mathbb{T})$,
\begin{equation}\label{equ 4.9}
\lim_{\delta\rightarrow 0}\limsup_{N\rightarrow+\infty}\sup_{\upsilon\in \mathcal{T}, s\leq\delta}P\big(|V^N_{\upsilon+s}(f)-V^N_{\upsilon}(f)|>\epsilon\big)=0,
\end{equation}
where $\mathcal{T}$ is the set of stopping times of $\{X_t^N\}_{t\geq 0}$ bounded by $T$.

By Lemma \ref{lemma 3.1 approximate independence} and Equation \eqref{equ 4.2},
\[
\mathbb{E}\left(\left(V_t^N(f)\right)^2\right)=O(1)
\]
and then Equation \eqref{equ 4.8} follows from Markov's inequality.

We define
\[
\Upsilon_t^N=V_0^N(f)+\int_0^t\left(\mathcal{L}_N+\partial_s\right)V_s^N(f)ds\text{\quad and\quad}\Xi_t^N=V_t^N(f)-\Upsilon_t^N,
\]
then to prove Equation \eqref{equ 4.9} we only need to show that
\begin{equation}\label{equ 4.10}
\lim_{\delta\rightarrow 0}\limsup_{N\rightarrow+\infty}\sup_{\upsilon\in \mathcal{T}, s\leq\delta}P\big(|\Upsilon^N_{\upsilon+s}(f)-\Upsilon^N_{\upsilon}(f)|>\epsilon\big)=0
\end{equation}
and
\begin{equation}\label{equ 4.11}
\lim_{\delta\rightarrow 0}\limsup_{N\rightarrow+\infty}\sup_{\upsilon\in \mathcal{T}, s\leq\delta}P\big(|\Xi^N_{\upsilon+s}(f)-\Xi^N_{\upsilon}(f)|>\epsilon\big)=0.
\end{equation}
According to the definition of $\mathcal{L}_N$,
\begin{align}\label{equ 4.12}
\left(\mathcal{L}_N+\partial_s\right)V_s^N(f)=&\frac{1}{N^{\frac{3}{2}}}\sum_{i=1}^N\sum_{j\neq i}\sum_{k=0}^{+\infty}\left(X_s^N(i)-\mathbb{E}X_s^N(i)\right)\lambda_k\left(\frac{i}{N}, \frac{j}{N}\right)
\left(kf\left(\frac{j}{N}\right)-f\left(\frac{i}{N}\right)\right) \notag\\
&+\frac{1}{\sqrt{N}}\sum_{i=1}^N\sum_{k=0}^{+\infty}\left(X_s^N(i)-\mathbb{E}X_t^N(i)\right)\psi_k\left(\frac{i}{N}\right)(k-1)f\left(\frac{i}{N}\right).
\end{align}
Hence, by Equations \eqref{equ 4.2}, \eqref{equ 4.12} and Lemma \ref{lemma 3.1 approximate independence}, we have
\begin{equation}\label{equ 4.13}
\sup_{0\leq s\leq T+1}\mathbb{E}\left(\left(\left(\mathcal{L}_N+\partial_s\right)V_s^N(f)\right)^2\right)=O(1).
\end{equation}
Since
\[
\mathbb{E}\left(|\Upsilon^N_{\upsilon+s}(f)-\Upsilon^N_{\upsilon}(f)|^2\right)\leq \delta \int_0^{T+1}\mathbb{E}\left(\left((\mathcal{L}_N+\partial_s)V_s^N(f)\right)^2\right)ds
\]
for $\delta<1$ and $s<\delta$ according to Cauchy-Schwarz inequality, Equation \eqref{equ 4.10} follows from Markov's inequality and Equation \eqref{equ 4.13}.

According to Dynkin's martingale formula, $\{\Xi_t^N\}_{t\geq 0}$ is a martingale with quadratic variation process
\[
\langle \Xi^N\rangle_t=\int_0^t \left(\mathcal{L}_N\left(\left(V_s^N(f)\right)^2\right)-2V_s^N(f)\mathcal{L}_NV_s^N(f)\right)ds.
\]
According to the definition of $\mathcal{L}_N$,
\begin{align*}
&\mathcal{L}_N\left(\left(V_s^N(f)\right)^2\right)-2V_s^N(f)\mathcal{L}_NV_s^N(f)\\
&=\frac{1}{N^2}\sum_{i=1}^N\sum_{j\neq i}\sum_{k=0}^{+\infty}\lambda_k\left(\frac{i}{N}, \frac{j}{N}\right)X_s^N(i)\left(kf\left(\frac{j}{N}\right)-f\left(\frac{i}{N}\right)\right)^2\\
&\text{\quad}+\frac{1}{N}\sum_{i=1}^N\sum_{k=0}^{+\infty}X_s^N(i)\psi_k\left(\frac{i}{N}\right)(k-1)^2f^2\left(\frac{i}{N}\right).
\end{align*}
Therefore, by Equation \eqref{equ 4.2},
\begin{align}\label{equ 4.14}
\sup_{0\leq s\leq T+1}\mathbb{E}\left(\left(\mathcal{L}_N\left(\left(V_s^N(f)\right)^2\right)-2V_s^N(f)\mathcal{L}_NV_s^N(f)\right)^2\right)=O(1).
\end{align}
By Cauchy-Schwarz inequality,
\begin{align*}
&\mathbb{E}\left(\left|\Xi^N_{\upsilon+s}(f)-\Xi^N_{\upsilon}(f)\right|^2\right)
=\mathbb{E}\left(\langle\Xi^N\rangle_{\upsilon+s}-\langle\Xi^N\rangle_{\upsilon}\right)\\
&=\mathbb{E}\left(\int_\upsilon^{\upsilon+s}\left(\mathcal{L}_N\left(\left(V_r^N(f)\right)^2\right)-2V_r^N(f)\mathcal{L}_NV_r^N(f)\right)dr\right)\\
&\leq \sqrt{\delta\int_0^{T+1}\mathbb{E}\left(\left(\mathcal{L}_N\left(\left(V_r^N(f)\right)^2\right)-2V_r^N(f)\mathcal{L}_NV_r^N(f)\right)^2\right)ds}
\end{align*}
for $\delta<1$ and $s<\delta$. Therefore, Equation \eqref{equ 4.11} follows from Markov's inequality and Equation \eqref{equ 4.14}. Since Equations \eqref{equ 4.10} and \eqref{equ 4.11} hold, Equation \eqref{equ 4.9} holds and the proof is complete.

\qed

\proof[Proof of Lemma \ref{lemma 4.2}]

For each integer $k\geq 1$, let $J_k^N$ be the number of moments in $[0, t]$ when some particle dies and meanwhile gives birth to $k$ new particles, then to complete the proof we only need to show that
\begin{equation}\label{equ 4.5}
\lim_{N\rightarrow+\infty}P\left(\sum_{k\geq \epsilon N^{\frac{4+\varepsilon_0}{8+4\varepsilon_0}}}J_k^N\geq 1\right)=0.
\end{equation}
For $M>0$, conditioned on $\sup_{0\leq s\leq t}\sum_{i=1}^NX_s^N(i)\leq MN$, $\sum_{k\geq \epsilon N^{\frac{4+\varepsilon_0}{8+4\varepsilon_0}}}J_k^N$ is stochastic dominated from above by $Y\left(C_7^NMNt\right)$, where $\{Y(t)\}_{t\geq 0}$ is the Poisson process with rate $1$ and
\[
C_7^N=\left\|\sum_{k\geq N^{\frac{4+\varepsilon_0}{8+4\varepsilon_0}}\epsilon}\lambda_k\right\|_\infty+\left\|\sum_{k\geq N^{\frac{4+\varepsilon_0}{8+4\varepsilon_0}}\epsilon}\psi_k\right\|_\infty.
\]
Since
\begin{align*}
C_7^N&\leq \left\|\sum_{k\geq N^{\frac{4+\varepsilon_0}{8+4\varepsilon_0}}\epsilon}\frac{k^{2+\varepsilon_0}}{\left(N^{\frac{4+\varepsilon_0}{8+4\varepsilon_0}}\epsilon\right)^{2+\varepsilon_0}}\lambda_k\right\|_\infty+\left\|\sum_{k\geq N^{\frac{4+\varepsilon_0}{8+4\varepsilon_0}}\epsilon}\frac{k^{2+\varepsilon_0}}{\left(N^{\frac{4+\varepsilon_0}{8+4\varepsilon_0}}\epsilon\right)^{2+\varepsilon_0}}\psi_k\right\|_\infty\\
&\leq \frac{\left\|\sum_{k\geq 0}k^{2+\varepsilon_0}\lambda_k\right\|_\infty+\left\|\sum_{k\geq 0}k^{2+\varepsilon_0}\psi_k\right\|_\infty}{N^{1+\frac{\varepsilon_0}{4}}\epsilon^{2+\varepsilon_0}},
\end{align*}
by Markov's inequality we have
\begin{align}\label{equ 4.6}
&P\Bigg(\sum_{k\geq \epsilon\sqrt{N}}J_k^N\geq 1\Bigg|\sup_{0\leq s\leq t}\sum_{i=1}^NX_s^N(i)\leq MN\Bigg) \notag\\
&\leq \frac{\left\|\sum_{k\geq 0}k^{2+\varepsilon_0}\lambda_k\right\|_\infty+\left\|\sum_{k\geq 0}k^{2+\varepsilon_0}\psi_k\right\|_\infty}{N^{\frac{\varepsilon_0}{4}}\epsilon^{2+\varepsilon_0}}Mt.
\end{align}
Let $\hat{X}_t^N$ be our process with parameters $\{\hat{\lambda}_k\}_{k\geq 0}$ and $\{\hat{\psi}_k\}_{k\geq 0}$ such that $\hat{\lambda}_0=0, \hat{\psi}_0=0$ and $\hat{\lambda}_k=\lambda_k, \hat{\psi}_k=\psi_k$ for all $k\geq 1$, then $\sum_{i=1}^N\hat{X}_t^N(i)$ is increasing in $t$. Let $\hat{X}_0^N=X_0^N$, then $\sum_{i=1}^{N}X_t^N(i)$ is stochastic dominated from above by $\sum_{i=1}^N\hat{X}_t^N(i)$ and hence
\[
P\left(\sup_{0\leq s\leq t}\sum_{i=1}^NX_s^N(i)\geq MN\right)\leq P\left(\sum_{i=1}^N\hat{X}_t^N(i)\geq NM\right).
\]
Then by $\{\hat{X}_t^N\}_{t\geq 0}$-version Equation (4.1) and Markov's inequality,
\begin{equation}\label{equ 4.7}
\lim_{M\rightarrow+\infty}\limsup_{N\rightarrow+\infty}P\left(\sup_{0\leq s\leq t}\sum_{i=1}^NX_s^N(i)\geq MN\right)=0.
\end{equation}
Equation \eqref{equ 4.5} follows from Equations \eqref{equ 4.6} and \eqref{equ 4.7} and the proof is complete.

\qed

At last, we prove Theorem \ref{theorem 2.2 fluctuation}.

\proof[Proof of Theorem \ref{theorem 2.2 fluctuation}]

By Lemma \ref{lemma 4.1}, any subsequence of $\{V_t^N:~0\leq t\leq T\}_{N\geq 1}$ has weakly convergent subsequence. Let $\{\hat{V}_t\}_{0\leq t\leq T}$ be a weak limit of a subsequence of $\{V_t^N:~0\leq t\leq T\}_{N\geq 1}$, then we only need to show that $\{\hat{V}_t\}_{0\leq t\leq T}=\{V_t\}_{0\leq t\leq T}$ to complete the proof. In this proof we still denote the subsequence convergent to $\{\hat{V}_t\}_{0\leq t\leq T}$ by $\{V_t^N:~0\leq t\leq T\}_{N\geq 1}$ for simplicity. By Lemma \ref{lemma 4.2},
\[
\lim_{N\rightarrow+\infty}P\left(\sup_{0\leq t\leq T}\left|V_t^N(f)-V_{t-}^N(f)\right|\geq \frac{\epsilon}{N^{\frac{\varepsilon_0}{8+4\varepsilon_0}}}\right)=0
\]
for any $f\in C^{\infty}(\mathbb{T})$ and $\epsilon>0$. Hence, $\{\hat{V}_t(f)\}_{0\leq t\leq T}$ is continuous in $t$ for any $f\in C^\infty(\mathbb{T})$. According to Assumption (A), $\hat{V}_0(f)$ follows from $\mathbb{N}\left(0, \int_{\mathbb{T}}\phi(u)f^2(u)du\right)$. Therefore, to complete the proof we only need to show that $\{m_t(G,f)\}_{0\leq t\leq T}$ is a martingale for any $G\in C_c^{\infty}(\mathbb{R})$ and $f\in C^\infty(\mathbb{T})$, where
\begin{align*}
m_t(G,f)=G\left(\hat{V}_t(f)\right)-G\left(\hat{V}_0(f)\right)&-\int_0^tG^\prime\left(\hat{V}_s(f)\right)\hat{V}_s\left(\mathcal{P}_1f\right)ds\\
&-\frac{1}{2}\int_0^tG^{\prime\prime}\left(\hat{V}_s(f)\right)\left(\|\mathcal{A}_sf\|_2^2+\sum_{k=0}^{+\infty}\|\mathcal{U}_s^kf\|_2^2\right)ds.
\end{align*}
For each $N\geq 1$, let
\[
m_t^N(G,f)=G\left(V_t^N(f)\right)-G\left(V_0^N(f)\right)-\int_0^t(\mathcal{L}_N+\partial_s)G\left(V_s^N(f)\right)ds,
\]
then $\{m_t^N(G,f)\}_{0\leq t\leq T}$ is a martingale by Dynkin's martingale formula. According to the definition of $\mathcal{L}_N$, Equation \eqref{equ 4.2} and Taylor's expansion up to the second order with Lagrange's remainder,
\begin{align*}
m_t^N(G,f)=&G\left(V_t^N(f)\right)-G\left(V_0^N(f)\right) \\
&-\int_0^tG^\prime\left(V_s^N(f)\right){\rm \uppercase\expandafter{\romannumeral1}}_s^Nds
-\frac{1}{2}\int_0^tG^{\prime\prime}\left(V_s^N(f)\right){\rm \uppercase\expandafter{\romannumeral2}}_s^Nds+\delta_{4,t}^N
\end{align*}
for $0\leq t\leq T$, where
\[
\sup_{0\leq t\leq T} \mathbb{E}\left(\left(\delta^N_{4,t}\right)^2\right)\leq \frac{C_{10}\|G^{\prime\prime\prime}\|^2_\infty}{N}
\]
for some $C_{10}<+\infty$ independent of $N$,
\begin{align*}
{\rm \uppercase\expandafter{\romannumeral1}}_s^N=&\frac{1}{\sqrt{N}}\sum_{i=1}^N\left(X_s^N(i)-\mathbb{E}X_s^N(i)\right)
\frac{1}{N}\sum_{j\neq i}\sum_{k=0}^{+\infty}\lambda_k\left(\frac{i}{N}, \frac{j}{N}\right)\left(kf\left(\frac{j}{N}\right)-f\left(\frac{i}{N}\right)\right)\\
&+\frac{1}{\sqrt{N}}\sum_{i=1}^N\left(X_s^N(i)-\mathbb{E}X_s^N(i)\right)\sum_{k=0}^{+\infty}\psi_k\left(\frac{i}{N}\right)(k-1)f\left(\frac{i}{N}\right)
\end{align*}
and
\begin{align*}
{\rm \uppercase\expandafter{\romannumeral2}}_s^N=&\frac{1}{N}\sum_{i=1}^NX_s^N(i)\frac{1}{N}\sum_{j\neq i}\sum_{k=0}^{+\infty}\lambda_k\left(\frac{i}{N}, \frac{j}{N}\right)\left(kf\left(\frac{j}{N}\right)-f\left(\frac{i}{N}\right)\right)^2\\
&+\frac{1}{N}\sum_{i=1}^NX_s^N(i)\sum_{k=0}^{+\infty}\psi_k\left(\frac{i}{N}\right)(k-1)^2f^2\left(\frac{i}{N}\right).
\end{align*}
Therefore, $\{m_t^N(G, f)\}_{N\geq 1}$ are uniformly integrable for any $0\leq t\leq T$ according to Lemma \ref{lemma 3.1 approximate independence} and Equation \eqref{equ 4.2}. Hence, by Theorem 5.3 of \cite{Whitt2007}, to prove that $\{m_t(G,f)\}_{0\leq t\leq T}$ is a martingale we only need to show that $m_t^N(G,f)$ converges weakly to $m_t(G,f)$ as $N\rightarrow+\infty$ for any $0\leq t\leq T$. By the definition of ${\rm \uppercase\expandafter{\romannumeral1}}_s^N$, Lemma \ref{lemma 3.1 approximate independence} and Equation \eqref{equ 4.2},
\[
{\rm \uppercase\expandafter{\romannumeral1}}_s^N=V_s^N(\mathcal{P}_1^Nf)+\delta^N_{5,s}
\]
and
\[
{\rm \uppercase\expandafter{\romannumeral2}}_s^N=\mu_s^N(\mathcal{P}_3^Nf)+\delta^N_{6,s},
\]
where
\[
\sup_{0\leq s\leq T}\mathbb{E}\left(\left(\delta^N_{5,s}\right)^2\right)=O(N^{-1}), \text{\quad} \sup_{0\leq s\leq T}\mathbb{E}\left(\left(\delta^N_{6,s}\right)^2\right)
=O(N^{-2}),
\]
\begin{align*}
(\mathcal{P}_1^Nf)(u)=\frac{1}{N}\sum_{j=1}^N\sum_{k=0}^{+\infty}\lambda_k\left(u, \frac{j}{N}\right)\left(kf\left(\frac{j}{N}\right)-f(u)\right)+\sum_{k=0}^{+\infty}\psi_k(u)(k-1)f(u)
\end{align*}
and
\begin{align*}
(\mathcal{P}_3^Nf)(u)=\frac{1}{N}\sum_{j=1}^N\sum_{k=0}^{+\infty}\lambda_k\left(u, \frac{j}{N}\right)\left(kf\left(\frac{j}{N}\right)-f(u)\right)^2+\sum_{k=0}^{+\infty}\psi_k(u)(k-1)^2f^2(u)
\end{align*}
for any $u\in \mathbb{T}$. Since $f\in C^\infty(\mathbb{T})$,
\[
\|\mathcal{P}_lf-\mathcal{P}_l^Nf\|_\infty=O(N^{-1})
\]
for $l=1,3$, where
\[
(\mathcal{P}_3f)(u)=\int_{\mathbb{T}}\sum_{k=0}^{+\infty}\lambda_k(u,v)\left(kf(v)-f(u)\right)^2dv+\sum_{k=0}^{+\infty}\psi_k(u)(k-1)^2f^2(u)
\]
for any $u\in \mathbb{T}$. Therefore, by Theorem \ref{theorem 2.1 scaling limit}, Lemma \ref{lemma 3.1 approximate independence}
and Equation \eqref{equ 4.2},
\[
{\rm \uppercase\expandafter{\romannumeral1}}_s^N=V_s^N(\mathcal{P}_1f)+\delta^N_{7,s}
\]
and
\[
{\rm \uppercase\expandafter{\romannumeral2}}_s^N=\mu_s(\mathcal{P}_3f)+\delta^N_{8,s}
=\|\mathcal{A}_sf\|_2^2+\sum_{k=0}^{+\infty}\|\mathcal{U}_s^kf\|_2^2+\delta^N_{8,s},
\]
where
\[
\lim_{N\rightarrow+\infty}\int_0^T|\delta^N_{7,s}|ds=0 \text{\quad and\quad}\lim_{N\rightarrow+\infty}\int_0^T|\delta^N_{8,s}|ds=0
\]
in probability. In conclusion,
\begin{align*}
m_t^N(G,f)=&G\left(V_t^N(f)\right)-G\left(V_0^N(f)\right)-\int_0^tG^\prime\left(V_s^N(f)\right)V_s^N\left(\mathcal{P}_1f\right)ds\\
&-\frac{1}{2}\int_0^tG^{\prime\prime}\left(V_s^N(f)\right)\left(\|\mathcal{A}_sf\|_2^2+\sum_{k=0}^{+\infty}\|\mathcal{U}_s^kf\|_2^2\right)ds+\delta^N_{9,t},
\end{align*}
where $\lim_{N\rightarrow+\infty}\delta^N_{9,t}=0$ in probability. Let $N\rightarrow+\infty$ in the above equation, then $m_t^N(G,f)$ converges weakly to $m_t(G,f)$ and the proof is complete.

\qed

\section{Proof of Theorem \ref{theorem 2.3 hittingTimes}} \label{section five}
In this section we prove Theorem \ref{theorem 2.3 hittingTimes}. Throughout this section we assume that $f\in C^\infty(\mathbb{T})$ and $r>\int_{\mathbb{T}}\phi(u)du$ make $\tau_{r,f}<+\infty$ and $\mu_{\tau_{r,f}}(\mathcal{P}_1f)>0$.

\proof[Proof of Theorem \ref{theorem 2.3 hittingTimes}]

As we have introduced in Remark \ref{remark 2.3}, $\frac{d}{dt}\mu_t(f)=\mu_t(\mathcal{P}_1f)$. Since
\[
\mu_{\tau_{r,f}}(\mathcal{P}_1f)=\frac{d}{dt}\mu_t(f)\Big|_{t=\tau_{r,f}}>0,
\]
there exists $\delta_1>0$ such that $\frac{d}{dt}\mu_t(f)\Big|_{t=s}>0$ for $s\in \left[\tau_{r,f}-\delta_1, \tau_{r,f}+\delta_1\right]$ and then $\mu_t(f)$ is strictly increasing in $s\in\left[\tau_{r,f}-\delta_1, \tau_{r,f}+\delta_1\right]$. For any $\epsilon\in (0, \delta_1)$,
\[
\left\{\tau_{r,f}^N>\tau_{r,f}+\epsilon\right\}\subseteq \left\{\mu^N_{\tau_{r,f}+\epsilon}<r\right\}
\]
and hence
\[
P\left(\tau_{r,f}^N>\tau_{r,f}+\epsilon\right)\leq P\left(\mu^N_{\tau_{r,f}+\epsilon}<r\right).
\]
By Theorem \ref{theorem 2.1 scaling limit},
\[
\lim_{N\rightarrow+\infty}\mu^N_{\tau_{r,f}+\epsilon}=\mu_{\tau_{r,f}+\epsilon}>\mu_{\tau_{r,f}}=r
\]
in probability and hence
\begin{equation}\label{equ 5.1}
\lim_{N\rightarrow+\infty}P\left(\tau_{r,f}^N>\tau_{r,f}+\epsilon\right)=0.
\end{equation}
For any $T>0$, consider $\{\mu_t^N(f)\}_{0\leq t\leq T}$ as a random element in $\mathcal{D}\left([0, T], \mathbb{R}\right)$, then we claim that $\{\mu_t^N(f):~0\leq t\leq T\}_{N\geq 1}$ are tight. We prove this claim at the end of this section. Then, by Theorem \ref{theorem 2.1 scaling limit}, $\{\mu_t^N(f)\}_{0\leq t\leq T}$ converges weakly to $\{\mu_t(f)\}_{0\leq t\leq T}$ as $N\rightarrow+\infty$ and hence
\[
\lim_{N\rightarrow+\infty}\sup_{0\leq t\leq T}\mu_t^N(f)=\sup_{0\leq t\leq T}\mu_t(f)
\]
in probability. Since $\left\{\tau_{r,f}^N<\tau_{r,f}-\epsilon\right\}\subseteq \left\{\sup_{0\leq t\leq \tau_{r,f}-\epsilon}\mu_t^N(f)\geq r\right\}$ and
\[
\sup_{0\leq t\leq \tau_{r,f}-\epsilon}\mu_t(f)<r,
\]
we have
\begin{equation}\label{equ 5.2}
\lim_{N\rightarrow+\infty}P\left(\tau_{r,f}^N<\tau_{r,f}-\epsilon\right)=0.
\end{equation}
By Equations \eqref{equ 5.1} and \eqref{equ 5.2}, $\lim_{N\rightarrow+\infty}\tau_{r,f}^N=\tau_{r,f}$ in probability.

\quad

For each $N\geq 1$, let $\delta_2^N=\sqrt{N}\left(\mu^N_{\tau^N_{r,f}}(f)-r\right)$, then
\[
\delta_2^N\leq \sqrt{N}\left(\mu^N_{\tau^N_{r,f}}(f)-\mu^N_{\tau^N_{r,f}-}(f)\right)
\]
and hence
\[
\left\{\delta_2^N\geq \epsilon\right\}\subseteq \left\{\tau^N_{r,f}>\tau_{r,f}+1\right\}\bigcup\left\{\sup_{0\leq s\leq \tau_{r,f}+1}\left|\mu_s^N(f)-\mu_{s-}^N(f)\right|>\frac{\epsilon}{\sqrt{N}}\right\}
\]
for any $\epsilon>0$. Therefore, by Equation \eqref{equ 5.1} and Lemma \ref{lemma 4.2},
\begin{equation}\label{equ 5.3}
\lim_{N\rightarrow+\infty}\delta_2^N=0
\end{equation}
in probability. Let $\delta_3^N=\left|\sqrt{N}\left(\mu_{\tau^N_{r,f}}(f)-\varpi_{\tau^N_{r,f}}(f)\right)\right|$. By Equation \eqref{equ 4.4},
\[
\sup_{0\leq t\leq T}\left|\mu_t(f)-\varpi_t^N(f)\right|=O(N^{-1})
\]
for given $f\in C^\infty(\mathbb{T})$ and $T>0$. Hence, for any $\epsilon>0$,
\[
\sup_{0\leq t\leq \tau_{r,f}+1}\left|\sqrt{N}\left(\mu_t(f)-\varpi_t^N(f)\right)\right|<\epsilon
\]
and
\[
\left\{\delta_3^N\geq \epsilon\right\}\subseteq \left\{\tau^N_{r,f}>\tau_{r,f}+1\right\}
\]
when $N$ is sufficiently large. Therefore, by Equation \eqref{equ 5.1},
\begin{equation}\label{equ 5.4}
\lim_{N\rightarrow+\infty}\delta_3^N=0
\end{equation}
in probability. Since $\mu_{\tau_{r,f}}(f)=r$, by Equations \eqref{equ 5.3} and \eqref{equ 5.4},
\begin{equation}\label{equ 5.5}
V_{\tau^N_{r,f}}^N(f)=-\sqrt{N}\left(\mu_{\tau^N_{r,f}}(f)-\mu_{\tau_{r,f}}(f)\right)+\delta_4^N
\end{equation}
with $\lim_{N\rightarrow+\infty}\delta_4^N=0$ in probability. By Lagrange's mean value theorem,
\begin{align*}
\mu_{\tau^N_{r,f}}(f)-\mu_{\tau_{r,f}}(f)&=\frac{d}{dt}\mu_t(f)\Big|_{t={\sigma_N}}\left(\tau^N_{r,f}-\tau_{r,f}\right)\\
&=\mu_{\sigma_N}(\mathcal{P}_1f)\left(\tau^N_{r,f}-\tau_{r,f}\right)
\end{align*}
for some $\sigma_N$ between $\tau^N_{r,f}$ and $\tau_{r,f}$. Since $\lim_{N\rightarrow+\infty}\tau^N_{r,f}=\tau_{r,f}$ in probability,
\[
\lim_{N\rightarrow+\infty}\mu_{\sigma_N}(\mathcal{P}_1f)=\mu_{\tau_{r,f}}(\mathcal{P}_1f)
\]
in probability. Therefore, by Corollary \ref{corollary 2.3} and Equation \eqref{equ 5.5}, to complete the proof we only need to show that $V^N_{\tau^N_{r,f}}(f)$ converges weakly to $V_{\tau_{r,f}}(f)$ as $N\rightarrow+\infty$. For any bounded continuous function $H$ and $\epsilon>0$,
\begin{align*}
\mathbb{E}H(V^N_{\tau^N_{r,f}}(f))&=\mathbb{E}\left(H(V^N_{\tau^N_{r,f}}(f))1_{\{|\tau_{r,f}^N-\tau_{r,f}|\leq \epsilon\}}\right)+o(1)\\
&\leq \mathbb{E}\left(\sup_{\tau_{r,f}-\epsilon\leq s\leq \tau_{r,f}+\epsilon}H(V^N_s(f))\right)+o(1).
\end{align*}
Then, by Theorem \ref{theorem 2.2 fluctuation},
\[
\limsup_{N\rightarrow+\infty}\mathbb{E}H(V^N_{\tau^N_{r,f}}(f))\leq \mathbb{E}\left(\sup_{\tau_{r,f}-\epsilon\leq s\leq \tau_{r,f}+\epsilon}H(V_s(f))\right).
\]
Since $\epsilon$ is arbitrary and $V_t(f)$ is continuous in $t$, let $\epsilon\rightarrow 0$, then
\[
\limsup_{N\rightarrow+\infty}\mathbb{E}H(V^N_{\tau^N_{r,f}}(f))\leq \mathbb{E}H(V_{\tau_{r,f}}(f)).
\]
Since $H$ is arbitrary, replace $H$ by $-H$ in the above inequality, then
\[
\liminf_{N\rightarrow+\infty}\mathbb{E}H(V^N_{\tau^N_{r,f}}(f))\geq \mathbb{E}H(V_{\tau_{r,f}}(f))
\]
and hence
\[
\lim_{N\rightarrow+\infty}\mathbb{E}H(V^N_{\tau^N_{r,f}}(f))=\mathbb{E}H(V_{\tau_{r,f}}(f)).
\]

As a result, $V^N_{\tau^N_{r,f}}(f)$ converges weakly to $V_{\tau_{r,f}}(f)$ as $N\rightarrow+\infty$ and the proof is complete.

\qed

At last, we prove the tightness of $\{\mu_t^N(f):~0\leq t\leq T\}_{N\geq 1}$.

\proof[Proof of the tightness of $\{\mu_t^N(f):~0\leq t\leq T\}_{N\geq 1}$]

By Aldous' criterion, we only need to check following two Equations.

(1) For any $t>0$,
\begin{equation}\label{equ 5.6}
\lim_{M\rightarrow+\infty}\limsup_{N\rightarrow+\infty}P\big(|\mu_t^N(f)|\geq M\big)=0.
\end{equation}

(2) For any $\epsilon>0$,
\begin{equation}\label{equ 5.7}
\lim_{\delta\rightarrow 0}\limsup_{N\rightarrow+\infty}\sup_{\upsilon\in \mathcal{T}, s\leq\delta}P\big(|\mu^N_{\upsilon+s}(f)-\mu^N_{\upsilon}(f)|>\epsilon\big)=0,
\end{equation}
where $\mathcal{T}$ is the set of stopping times of $\{X_t^N\}_{t\geq 0}$ bounded by $T$ as we have introduced in Section \ref{section four}.

By Markov's inequality, $P\big(|\mu_t^N(f)|\geq M\big)\leq \frac{\mathbb{E}\left|\mu_t^N(f)\right|}{M}$ and hence Equation \eqref{equ 5.6} follows from Equation \eqref{equ 4.1}. To check Equation \eqref{equ 5.7}, we define
\[
Z_t^N=\mu_0^N(f)+\int_0^t\mathcal{L}_N\mu_s^N(f)ds\text{\quad and\quad}\Omega_t^N=\mu_t^N(f)-Z_t^N,
\]
then Equation \eqref{equ 5.7} follows from
\begin{equation}\label{equ 5.8}
\lim_{\delta\rightarrow 0}\limsup_{N\rightarrow+\infty}\sup_{\upsilon\in \mathcal{T}, s\leq\delta}P\big(|Z^N_{\upsilon+s}(f)-Z^N_{\upsilon}(f)|>\epsilon\big)=0
\end{equation}
and
\begin{equation}\label{equ 5.9}
\lim_{\delta\rightarrow 0}\limsup_{N\rightarrow+\infty}\sup_{\upsilon\in \mathcal{T}, s\leq\delta}P\big(|\Omega^N_{\upsilon+s}(f)-\Omega^N_{\upsilon}(f)|>\epsilon\big)=0.
\end{equation}
By direct calculation
\begin{align*}
\mathcal{L}_N\mu_s^N(f)=&\frac{1}{N^2}\sum_{i=1}^N\sum_{j\neq i}\sum_{k=0}^{+\infty}X_s^N(i)\lambda_k\left(\frac{i}{N}, \frac{j}{N}\right)
\left(kf\left(\frac{j}{N}\right)-f\left(\frac{i}{N}\right)\right)\\
&+\frac{1}{N}\sum_{i=1}^N\sum_{k=0}^{+\infty}X_s^N(i)\psi_k\left(\frac{i}{N}\right)(k-1)f\left(\frac{i}{N}\right).
\end{align*}
Then, by Equation \eqref{equ 4.2}, there exists $C_9<+\infty$ independent of $N$ such that
\begin{equation}\label{equ 5.10}
\mathbb{E}\left(\left(\mathcal{L}_N\mu_s^N(f)\right)^2\right)\leq C_9\|f\|_\infty^2
\end{equation}
for all $0\leq s\leq T+1$. By Cauchy-Schwarz inequality,
\[
\mathbb{E}\left(|Z^N_{\upsilon+s}(f)-Z^N_{\upsilon}(f)|^2\right)\leq \delta \int_0^{T+1}\mathbb{E}\left(\left(\mathcal{L}_N\mu_s^N(f)\right)^2\right)ds
\]
for $\delta<1$ and $s<\delta$, then Equation \eqref{equ 5.8} follows from Markov's inequality and Equation \eqref{equ 5.10}.

According to Dynkin's martingale formula, $\{\Omega_t^N\}_{t\geq 0}$ is a martingale with quadratic variation process
\[
\langle \Omega^N\rangle_t=\int_0^t \left(\mathcal{L}_N\left(\left(\mu_s^N(f)\right)^2\right)-2\mu_s^N(f)\mathcal{L}_N\mu_s^N(f)\right)ds.
\]
By direct calculation,
\begin{align*}
&\mathcal{L}_N\left(\left(\mu_s^N(f)\right)^2\right)-2\mu_s^N(f)\mathcal{L}_N\mu_s^N(f)\\
&=\frac{1}{N^3}\sum_{i=1}^N\sum_{j\neq i}\sum_{k=0}^{+\infty}\lambda_k\left(\frac{i}{N}, \frac{j}{N}\right)X_s^N(i)\left(kf\left(\frac{j}{N}\right)-f\left(\frac{i}{N}\right)\right)^2\\
&\text{\quad}+\frac{1}{N^2}\sum_{i=1}^N\sum_{k=0}^{+\infty}X_s^N(i)\psi_k\left(\frac{i}{N}\right)(k-1)^2f^2\left(\frac{i}{N}\right).
\end{align*}
Therefore, by Equation \eqref{equ 4.1},
\begin{align}\label{equ 5.11}
&\mathbb{E}\left(\left|\Omega^N_{\upsilon+s}-\Omega^N_\upsilon\right|^2\right)=\mathbb{E}\left(\langle \Omega^N\rangle_{\upsilon+s}-\langle \Omega^N\rangle_\upsilon\right)\\
&\leq \int_0^{T+1}\mathbb{E}\left(\mathcal{L}_N\left(\left(\mu_s^N(f)\right)^2\right)-2\mu_s^N(f)\mathcal{L}_N\mu_s^N(f)\right)ds=O(N^{-1}) \notag
\end{align}
for $\delta<1$ and $s<\delta$. As a result, Equation \eqref{equ 5.9} follows from Markov's inequality and Equation \eqref{equ 5.11}. Since Equations \eqref{equ 5.8} and \eqref{equ 5.9} hold, Equation \eqref{equ 5.7} holds and the proof is complete.

\qed

\section{Revisit to a density-dependent Markov chain case}\label{section six}
As an application of our main result, in this section we revisit a special case of our model which reduces to a density-dependent Markov chain. Throughout this section we assume that $\lambda_k\equiv 0$ for $k\neq 1$ and $\psi_k\equiv q_k$ for some $q_k\in [0,+\infty)$ for all $k\geq 0$. Let
\[
\mathcal{N}_t^N=\sum_{i=1}^NX_t^N(i),
\]
then $\{\mathcal{N}_t^N\}_{t\geq 0}$ is a density-dependent Markov chain as we have introduced in Section \ref{section one}. According to limit theorems of density-dependent Markov chains given in \cite{Kurtz1978}, we have following propositions.

\begin{proposition}\label{proposition 6.1}(Kurtz, 1978)
Under Assumption (A),
\[
\lim_{N\rightarrow+\infty}\frac{\mathcal{N}_t^N}{N}=n_t
\]
in probability for any $t\geq 0$, where
\[
\begin{cases}
&\frac{d}{dt}n_t=n_t\sum_{k=0}^{+\infty}(k-1)q_k,\\
&n_0=\int_{\mathbb{T}}\phi(u)du.
\end{cases}
\]
\end{proposition}

\begin{proposition}\label{proposition 6.2}(Kurtz, 1978)
Let $\alpha_t^N=\frac{\mathcal{N}_t^N-\mathbb{E}\mathcal{N}_t^N}{\sqrt{N}}$, then $\{\alpha_t^N\}_{0\leq t\leq T}$ converges weakly to $\{\alpha_t\}_{0\leq t\leq T}$ as $N\rightarrow+\infty$ for any $T\geq 0$, where
\[
\begin{cases}
&d\alpha_t=\left(\sum_{k=0}^{+\infty}(k-1)q_k\right)\alpha_tdt+\sqrt{\left(\sum_{k=0}^{+\infty}(k-1)^2q_k\right)n_t}dB_t, \\
&\alpha_0\text{~follows~}\mathbb{N}\left(0, \int_{\mathbb{T}}\phi(u)du\right),
\end{cases}
\]
where $\{B_t\}_{t\geq 0}$ is a standard Brownian motion.
\end{proposition}

We further assume that $\sum_{k=0}^{+\infty}(k-1)q_k>0$, then for any $r>\int_{\mathbb{T}}\phi(u)du$ we define $\beta_r=\inf\{t\geq 0:~n_t=r\}$ and $\beta_r^N=\inf\left\{t\geq 0:~\frac{\mathcal{N}_t^N}{N}\geq r\right\}$. By Theorem 11.4.1 of \cite{Ethier1986}, we have the following proposition.

\begin{proposition}\label{proposition 6.3}(Ethier and Kurtz, 1986)
For $r>\int_{\mathbb{T}}\phi(u)du$, $\sqrt{N}\left(\beta_r^N-\beta_r\right)$ converges weakly to $-\frac{\alpha_{\beta_r}}{\left(\sum_{k=0}^{+\infty}(k-1)q_k\right)n_{\beta_r}}$
as $N\rightarrow+\infty$.
\end{proposition}

Now we apply Theorems \ref{theorem 2.1 scaling limit}, \ref{theorem 2.2 fluctuation} and \ref{theorem 2.3 hittingTimes} to give alternative proofs of above propositions.

\proof[An alternative proof of Proposition \ref{proposition 6.1}]

Let $\vec{1}$ be the constant function taking value $1$, then $\frac{\mathcal{N}_t^N}{N}=\mu_t^N(\vec{1})$. Hence, by Theorem \ref{theorem 2.1 scaling limit},
\[
\lim_{N\rightarrow+\infty}\frac{\mathcal{N}_t^N}{N}=\int_{\mathbb{T}}\rho_t(u)\vec{1}(u)du=\int_{\mathbb{T}}\rho_t(u)du
\]
in probability, where $\rho_0=\phi$ and
\[
\frac{d}{dt}\rho_t=\mathcal{P}_2\rho_t
\]
such that
\begin{align*}
\int_{\mathbb{T}}\mathcal{P}_2\rho_t(u)du=&\left(\sum_{k=0}^{+\infty}(k-1)q_k\right)\int_{\mathbb{T}}\rho_t(u)du\\
&+\int_{\mathbb{T}^2}\lambda_1(v,u)\rho_t(v)dudv-\int_{\mathbb{T}^2}\lambda_1(u,v)\rho_t(u)dudv\\
=&\left(\sum_{k=0}^{+\infty}(k-1)q_k\right)\int_{\mathbb{T}}\rho_t(u)du.
\end{align*}
Therefore, $\int_{\mathbb{T}}\rho_t(u)du=n_t$ and the proof is complete.

\qed

\proof[An alternative proof of Proposition \ref{proposition 6.2}]

Since $\alpha_t^N=V_t^N(\vec{1})$, by Theorem \ref{theorem 2.2 fluctuation}, $\{\alpha_t^N\}_{0\leq t\leq T}$ converges weakly to $\{V_t(\vec{1})\}_{0\leq t\leq T}$ as $N\rightarrow+\infty$, where $V_0(\vec{1})$ follows $\mathbb{N}\left(0, \int_{\mathbb{T}}\phi(u)du\right)$ and
\[
dV_t(\vec{1})=V_t(\mathcal{P}_1(\vec{1}))dt+d\xi_t+d\eta_t
\]
with
\[
\xi_t=\left(\int_0^t\mathcal{A}_s^*d\mathcal{B}_s\right)(\vec{1}) \text{~and~}
\eta_t=\left(\int_0^t(\mathcal{U}_s^1)^*d\mathcal{W}_s^1\right)(\vec{1}).
\]
According to the definition of $\mathcal{U}_s^1$, $\mathcal{U}_s^1\vec{1}=0$ and hence $\eta_t=0$. According to the definition of $\mathcal{P}_1$, $\mathcal{P}_1\vec{1}=\left(\sum_{k=0}^{+\infty}(k-1)q_k\right)\vec{1}$ and
\[
V_t(\mathcal{P}_1(\vec{1}))=\left(\sum_{k=0}^{+\infty}(k-1)q_k\right)V_t(\vec{1}).
\]
According to the definition of $\mathcal{A}_s$,
\[
{\rm Cov}(\xi_t, \xi_t)=\int_0^t\|\mathcal{A}_s\vec{1}\|_2^2ds
\]
with
\[
\|\mathcal{A}_s\vec{1}\|_2^2=\sum_{k=0}^{+\infty}(k-1)^2q_k\int_{\mathbb{T}}\rho_s(u)du=\sum_{k=0}^{+\infty}(k-1)^2q_kn_s.
\]
Therefore, $d\xi_t=\sqrt{\left(\sum_{k=0}^{+\infty}(k-1)^2q_k\right)n_t}dB_t$ and
\[
dV_t(\vec{1})=\left(\sum_{k=0}^{+\infty}(k-1)q_k\right)V_t(\vec{1})dt+\sqrt{\left(\sum_{k=0}^{+\infty}(k-1)^2q_k\right)n_t}dB_t.
\]
As a result, $V_t(\vec{1})=\alpha_t$ and the proof is complete.

\qed

\proof[An alternative proof of Proposition \ref{proposition 6.3}]

As we have shown in the proof of Proposition \ref{proposition 6.1}, $n_t=\int_\mathbb{T}\rho_t(u)du=\mu_t(\vec{1})$. Then $\beta_r=\tau_{r, \vec{1}}$ and $\beta_r^N=\tau^N_{r, \vec{1}}$. Hence, by Theorem \ref{theorem 2.3 hittingTimes}, $\sqrt{N}\left(\beta_r^N-\beta_r\right)$ converges weakly to $\mathbb{N}\left(0, \frac{\theta^2_{\tau_{r, \vec{1}}}(\vec{1})}{\mu^2_{\tau_{r, \vec{1}}}\left(\mathcal{P}_1\vec{1}\right)}\right)$ as $N\rightarrow+\infty$. As we have shown in the proof of proposition \ref{proposition 6.2}, $\mathcal{P}_1\vec{1}=\left(\sum_{k=0}^{+\infty}(k-1)q_k\right)\vec{1}$ and hence
\[
\mu_{\tau_{r, \vec{1}}}\left(\mathcal{P}_1\vec{1}\right)=\left(\sum_{k=0}^{+\infty}(k-1)q_k\right)\mu_{\tau_{r, \vec{1}}}(\vec{1})=\left(\sum_{k=0}^{+\infty}(k-1)q_k\right)r=\left(\sum_{k=0}^{+\infty}(k-1)q_k\right)n_{\beta_r}.
\]
As a result, to complete the proof we only need to show that $\alpha_{\beta_r}$ follows $\mathbb{N}\left(0, \theta^2_{\tau_{r, \vec{1}}}(\vec{1})\right)$. According to definitions of $\{\alpha_t\}_{t\geq 0}$ and $n_t$, $n_t=n_0e^{t\left(\sum_{k=0}^{+\infty}(k-1)q_k\right)}=n_0e^{t\mathcal{P}_1}\vec{1}$ and
\[
\alpha_t=e^{t\left(\sum_{k=0}^{+\infty}(k-1)q_k\right)}\alpha_0
+\int_0^te^{(t-s)\left(\sum_{k=0}^{+\infty}(k-1)q_k\right)}\sqrt{\left(\sum_{k=0}^{+\infty}(k-1)^2q_k\right)n_s}dB_s.
\]
Hence, $\alpha_{\beta_r}$ follows normal distribution with mean $0$ and variance
\[
{\rm Var}(\alpha_{\beta_r})=\frac{n^2_{\beta_r}}{n_0^2}\int_{\mathbb{T}}\phi(u)du
+\frac{1}{n_0^2}\int_0^{\beta_r}n_{2\beta_r-s}\left(\sum_{k=0}^{+\infty}(k-1)^2q_k\right)ds.
\]
As we have shown in the proof of Proposition \ref{proposition 6.2},
\[
\mathcal{U}_s^1e^{(t-s)\mathcal{P}_1}\vec{1}=e^{(t-s)\sum_{k=0}^{+\infty}(k-1)q_k}\mathcal{U}_s^1=0
\]
and
\begin{align*}
\|\mathcal{A}_se^{(t-s)\mathcal{P}_1}\vec{1}\|_2^2&=(e^{(t-s)\sum_{k=0}^{+\infty}(k-1)q_k})^2\|\mathcal{A}_s\vec{1}\|_2^2\\
&=\frac{1}{n_0^2}n^2_{t-s}n_s\sum_{k=0}^{+\infty}(k-1)^2q_k=\frac{1}{n_0^2}n_{2t-s}\sum_{k=0}^{+\infty}(k-1)^2q_k.
\end{align*}
Therefore,
\[
\theta^2_{\tau_{r, \vec{1}}}(\vec{1})=\frac{n^2_{\tau_{r, \vec{1}}}}{n_0^2}\int_{\mathbb{T}}\phi(u)du
+\frac{1}{n_0^2}\int_0^{\tau_{r, \vec{1}}}n_{2\tau_{r,\vec{1}}-s}\sum_{k=0}^{+\infty}(k-1)^2q_kds.
\]
As a result, ${\rm Var}(\alpha_{\beta_r})=\theta^2_{\tau_{r, \vec{1}}}(\vec{1})$ follows from $\tau_{r, \vec{1}}=\beta_r$ and the proof is complete.

\qed

\quad

\textbf{Data Availability.} The authors declare that all data supporting this article are available within the article.

\quad

\textbf{Acknowledgments.} The authors are grateful to Prof. Hui He for useful suggestions and comments. The authors are grateful to the financial
support from Beijing Jiaotong University with grant number 2022JBMC039.

{}
\end{document}